\documentclass[report]{svmult}
\makeatletter
\providecommand{\@LN}[2]{}
\makeatother
\usepackage[toc,page]{appendix}
\usepackage[utf8x]{inputenc}
\usepackage[T1]{fontenc}
\usepackage{algorithm}
\usepackage[noend]{algpseudocode}
\usepackage{url}
\usepackage[square,sort,comma,numbers]{natbib}
\usepackage{amsmath}
\usepackage{graphicx}
\usepackage[colorinlistoftodos]{todonotes}
\usepackage{subcaption}
\usepackage{caption}

\usepackage{bm}
\usepackage{amssymb}
\usepackage{amsfonts}
\DeclareMathAlphabet{\mathitbf}{OML}{cmm}{b}{it}
\def\jmath{j}

\newcommand{\gPC}{gPC\;}
\newcommand{\QMC}{qMC\;}
\newcommand{\mysol}{c}

\newcommand{\bsol}{\mathbf{\mysol}}
\newcommand{\Do}{\mathcal{D}}
\newcommand{\myug}{ug4\;}
\newcommand{\refeq}[1]{(\ref{#1})}

\newcommand{\bzero}{\textbf{0}}
\newcommand{\vTheta}{\varTheta}

\newcommand{\Pol}{\Psi}
\newcommand{\pol}{\psi}
\newcommand{\EXP}[1]{\mathbb{E}\left(#1\right)}
\newcommand{\di}{\mathrm{d}}
\newcommand{\btheta}{\mbox{\boldmath$\theta$}}

\def\sol{c}
\def\bsol{\mathbf{\sol}}
\newcommand{\Psii}{\bsol}

\def\D{\mathcal{D}}

\newcommand{\var}[1]{{\ensuremath{\mathrm{Var}}\mspace{-2mu}\left[#1\right]}}

\def\I{\mathbb{I}}

\newcommand{\err}{\mbox{err}}

\def\xib{\bm{\xi}}
\def\thetab{\bm{\theta}}

\def\bx{\mathbf{x}}

\def\be{\mathbf{e}}

\newcommand{\conc}{c} % mass fraction
\newcommand{\pres}{p} % hydrostatic pressure
\newcommand{\poro}{\phi} % porosity
\newcommand{\perm}{{\mathbf{K}}} % permeability
\newcommand{\dens}{\rho} % density
\newcommand{\visc}{\mu} % viscosity
\newcommand{\dvel}{{\mathbf{q}}} % Darcy velocity
\newcommand{\grav}{{\mathbf{g}}} % gravity
\newcommand{\disp}{{\mathbf{D}}} % diffusion-disperion
\begin{document}
%original was like this \graphicspath{{figures/}} % figures in subdirectory
%\graphicspath{{../figures/}} % figures in subdirectory
%\input{litvinenko} 
\title*{Solution of the 3D density-driven groundwater flow problem with uncertain porosity and permeability}
\author{Alexander Litvinenko\inst{1} \and Dmitry Logashenko\inst{2} \and Raul Tempone\inst{1,2} \and
Gabriel Wittum\inst{2,3} \and David Keyes\inst{2}}
 \authorrunning{Litvinenko et al.}
\institute{\inst{1}RWTH Aachen, Kackertstr. 9C, Aachen, Germany, Phone: +492418099203, \email{litvinenko, tempone@uq.rwth-aachen.de} \\ \inst{2}KAUST, Thuwal/Jeddah, Saudi Arabia, %Phone: +966128080685,
\email{dmitry.logashenko, gabriel.wittum, david.keyes@kaust.edu.sa} \\ \inst{3}G-CSC, Frankfurt University, Kettenhofweg 139, Frankfurt, Germany}
%\email{dmitry.logashenko, david.keyes, gabriel.wittum, raul.tempone@kaust.edu.sa}}
%
%
\maketitle

%THIS IS ONLINE ABSTRACT ONLY 
\abstract{
As groundwater is an essential nutrition and irrigation resource, its pollution may
lead to catastrophic consequences. Therefore, accurate modeling of the pollution of the soil and groundwater aquifer is highly important. As a model, we consider a density-driven groundwater flow problem with uncertain porosity and permeability. This problem may arise in geothermal reservoir simulation, natural saline-disposal basins, modeling of contaminant plumes, and subsurface flow. 
This strongly nonlinear time-dependent problem describes the convection of the two-phase flow. This liquid streams under the gravity force, building so-called ``fingers''. The accurate numerical solution requires fine spatial resolution with an unstructured mesh and, therefore, high computational resources. Consequently, we run the parallelized simulation toolbox \myug with the geometric multigrid solver on Shaheen II supercomputer.
The parallelization is done in physical and stochastic spaces. Additionally, we demonstrate how the \myug toolbox can be run in a black-box fashion for testing different scenarios in the density-driven flow.  
As a benchmark, we solve the Elder-like problem in a 3D domain. For approximations in the stochastic space, we use the generalized polynomial chaos expansion. We compute the mean, variance, and exceedance probabilities of the mass fraction. As a reference solution, we use the solution, obtained from the quasi-Monte Carlo method.}


%
%
%\begin{AMS}
%60H15,  %  Stochastic partial differential equations
%60H35,  %  Computational methods for stochastic equations 
%\end{AMS}

\section{Motivation}
\label{sec:1}

This work is a 3D extension of our previous 2D work \cite{LitvElder2D}. The problem setup is the same; the difference is only in the aquifer. In this work, we consider two 3D aquifers. These new 3D geometries contain a much larger number of degrees of freedom (DoFs), and, therefore, numerical experiments require significantly more computational resources.

Accurate prediction of the contamination of the groundwater is highly essential. Certain pollutants are soluble in water and can leak into groundwater systems, such as seawater into
coastal aquifers or wastewater leaks. Indeed, some pollutants can change the density of a fluid and induce density-driven flows
within the aquifer. This causes faster propagation of the contamination due to convection. Thus, simulation
and analysis of this density-driven flow play an important role in predicting how pollution can migrate through an aquifer \cite{Kobus2000, Density-driven_Fan97}.

In contrast to the transport of pollution by molecular diffusion, convection due to density-driven flow is
an unstable process that can undergo quite complicated patterns of distribution. 
The Elder problem is a simplified but comprehensive model that
describes the intrusion of salt water from a top boundary into an aquifer \cite{Elder_1967a, Elder_1967b,Voss_Souza,Elder_Overview17}. The evolution
of the concentration profile is typically referred to as \textit{fingering}. Note that due
to the nonlinear nature of the model, the distribution of the contamination strongly
depends on the model parameters, so that the system may have several stationary solutions
\cite{Johannsen2003}.

Hydrogeological formations typically have a complicated and heterogeneous structure, and geological media may consist of layers of porous media
with different porosities and permeability coefficients \cite{ScheiderKroehnPueschel2012, ReiterLogashenkoVogelWittum2017}. Difficulty in specifying hydrogeological parameters of the media, as
well as measuring the position and configuration of the layers, gives rise to certain errors. Typically, the averaged values of these
quantities are used. However, due to the nonlinearities of
the flow model, the averaging of the parameters does not necessarily lead to
the correct mathematical solution. Uncertainties arise from various factors, such as inaccurate measurements of the different parameters and the inability to
measure parameters at each spatial location at any given time. To model the uncertainties, we can use random variables, random fields, and random processes.
However, uncertainties in the input data spread through the model and make the solution
uncertain too. Thus, an accurate estimation of the output uncertainties is crucial.

Many techniques are used to investigate the propagation of uncertainties associated with porosity and permeability into the
solution, such as classical Monte Carlo sampling.
Although Monte Carlo sampling is dimension-independent, it converges slowly and therefore requires a lot of samples, and these simulations may be time-consuming and costly. Alternative, less costly techniques use
collocations, sparse grids, and surrogate methods, each with their advantages and disadvantages. But even advanced
techniques may require hundreds to thousands of simulations and assume a certain smoothness in the quantity of interest. Our simulations
contain up to 4.5MI spatial mesh points and 1000-3000 time-steps and are therefore run on a massive parallel cluster. Here, we
develop methods that require much fewer simulations, from a few dozen to a few hundred.

Perturbation methods, which decompose the quantity of interest (QoI) with respect to random parameters in a Taylor series, were considered and compared in  \cite{CREMER15_Fingers}.

Here, we use the well-established generalized polynomial chaos expansion (\gPC) technique \cite{xiuKarniadakis02a}, where we compute the \gPC coefficients by applying a Clenshaw-Curtis quadrature rule \cite{Barthelmann2000}. We use \gPC to compute QoIs such as the mean, the variance, and the exceedance probabilities of the mass fraction (in this case, a saltwater solution). We validate our obtained results using the quasi-Monte Carlo approach. Both methods require the computation of multiple simulations (scenarios) for variable porosity and permeability coefficients.

To the best of our knowledge, no other reported works have solved the Elder’s problem \cite{Voss_Souza,Elder_Overview17} with uncertain porosity and permeability parameters using \gPC.

Fingers in a highly unstable free convective flow have been studied in \cite{Xie12_Fingers}. 
Overviews of the
uncertainties in modeling groundwater solute transport \cite{OverviewUncert93} and modeling soil processes \cite{SoilOverview16} have been performed, as well as
reconnecting stochastic methods with hydrogeological applications \cite{NowakStochMethods18}, which included recommendations for optimization
and risk assessment. Fundamentals of stochastic hydrogeology, an overview of stochastic tools and accounting for uncertainty
have been described in \cite{rubin2003applHydro}.

Recent advances in uncertainty quantification, probabilistic risk assessment, and decision-making under uncertainty in hydrogeologic applications have been reviewed in \cite{TARTAKOVSKYI_Risk13}, where the
author reviewed probabilistic risk assessment methods in hydrogeology under parametric, geologic, and model uncertainties.
Density-driven vertical transport of saltwater through the freshwater lens has been modeled in \cite{POST17_Density-driven}.

Various methods can be used to compute the desired statistics, such as direct integration methods (Monte Carlo, quasi Monte Carlo, collocation)
and surrogate-based methods (generalized polynomial chaos approximation, stochastic Galerkin) \cite{Philipp12,matthiesKeese03cmame, babuska2004galerkin}.
Direct integration methods compute statistics by sampling 
unknown input coefficients (they are undefined and therefore uncertain)
and solving the corresponding PDE, while the surrogate-based method computes a low-cost functional (polynomial, exponential, trigonometrical) approximation of QoI.
Examples of the surrogate-based methods include radial basis functions \cite{liu2014,bompard2010,Loeven2007,giunta2004}, sparse polynomials \cite{
chkifa-adapt-stochfem-2015,Sudret_sparsePCE}, polynomial chaos expansion \cite{Xiu, Habib09_PCE, ConradMarzouk13, Dongbin} and low-rank tensor approximation \cite{DolgLitv15, Kim2006,dolgov2014computation,LitvSampling13,Litvin11PAMM}. A nonintrusive stochastic Galerkin was introduced in \cite{Giraldi2014}. The surrogate-based methods have an additional advantage if the gradients are available \cite{liu2017quantification}.  
Sparse grids methods to integrate high-dimensional integrals are considered in \cite{smoljak63, Griebel_Bungartz, Griebel, spiterp, novakRitter97, Novak1996, gerstnerGriebel98-numint, novakRitter99-simple, ConradMarzouk13, CONSTANTINE12}. The sparse grid software package is available in \cite{petrasSmolpak}.

The quantification of uncertainties in stochastic PDEs can pose a great challenge due to the possibly large number of random variables involved, and the high cost of each deterministic solution.

For problems with a large number of random variables, the methods based on a regular grid (in stochastic space) are less preferable due to the high computational cost. Instead, methods based on a scattered sampling scheme (MC, qMC) give more freedom in choosing the sample number $N$ and protect against sample failures. The MC quadrature and its variance-reduced variants have dimension-independent convergence $\mathcal{O}(N^{-\frac{1}{2}})$, and qMC has the worst case convergence $\mathcal{O}(\log^M(N)N^{-1})$, where $M$ is the dimension of the stochastic space \cite{matthies2007}. Collocations on sparse or full grids are affected by the dimension of the
integration domain \cite{babuska_collocation, nobile-sg-mc-2015, NobileTemponeWebster08}. 
In this work, we use the Halton rule to generate quasi-Monte Carlo quadrature points \cite{caf1998, joe2008}. A numerical comparison of other qMC sequences has been described in \cite{radovic1996}.

The construction of a \gPC-based surrogate model \cite{Matthies_encicl} is similar to computing the Fourier expansion, where only the Fourier coefficients need to be computed.
Some well-known functions, such as multivariate Legendre, Hermite, Chebyshev, and Laguerre functions, are taken as a basis \cite{xiuKarniadakis02a}. These functions should possess some useful properties, e.g., orthogonality. Usually, a
small number of quadrature points are used to compute these coefficients. Once the surrogate is constructed, its sampling could
be performed at a lower cost than a sampling of the original stochastic PDE.

To tackle the high numerical complexity, we implement a two-level parallelization.
First, we run all available quadrature or sampling points (also say scenarios) in parallel, with each scenario on a separate computing node. Second, each scenario is computed in parallel on $2-32$ cores.

This work is structured as follows: In Section \ref{litv:subsec:setup}, we outline the model of the density-driven groundwater flow in porous media and the numerical methods for this type of problem. We describe the stochastic modeling, integration methods, and the generalized polynomial chaos expansion technique in Section~\ref{litv:sec:StochModel}.
We present details of the parallelization of the computations in Section~\ref{sec:parallelization}.
Our multiple numerical results, described in Section \ref{litv:sec:numerics}, demonstrate the practical performance of the parallelized solver for the Elder-type problems with uncertainty coefficients in 3D computational domains. We conclude this work with a discussion in Section~\ref{litv:sec:Conclusion}.

\subsection{Our contribution} 
We applied the \gPC expansion to approximate the solution of the density driven flow. This is a time-dependent, nonlinear, second order differential equation. We estimated propagation of input uncertainties in the porosity and permeability into the mass fraction.

%
%
%%%%%%%%%%%%%%%%%%%%%%%%%%%%%%
\section{Problem setup}
\label{litv:subsec:setup}
\subsection{Density-driven groundwater flow}

The groundwater in the porous medium of the aquifers is subjected to gravitation. For this, any spatial variation of its density induces a flow. In this work, we regard the dependence of the density on the mass fraction of the dissolved salt. Traditional approaches for modelling of this system are described in \cite {Bear-2,Bear1979,Bear-Bachmat,DierschKolditz2002,POST17_Density-driven}. In this section, we briefly review the model to introduce the notation.

We consider a domain $\Do \subset \mathbb{R}^d$, $d=3$ consisting of two phases: the immobile solid porous matrix and a mobile liquid in its pores. We denote the porosity of the matrix by $\poro: \Do \to \mathbb{R}$ and its permeability $\perm: \Do \to \mathbb{R}^{d \times d}$ and the permeability by $\perm$. The liquid is a solution of $NaCl$ in water. The mass fraction of the brine (a saturated saltwater solution) in the liquid phase is $\conc (t, \mathbf{x}): [0, +\infty) \times \Do \to [0, 1]$. The density liquid phase is $\dens = \dens (\conc)$ and its viscosity $\visc = \visc (\conc)$.

The motion of the liquid phase through the solid matrix is characterized by the Darcy velocity $\dvel (t, \mathbf{x}): [0, +\infty) \times \Do \to \mathbb{R}^d$. Mass conservation of the liquid phase and the salt can be written as
\begin{eqnarray}
 \label {e_cont_eq}
 \partial_t (\poro \dens) & + & \nabla \cdot (\dens \dvel) = 0, \\
 \label {e_tran_eq}
 \partial_t (\poro \dens \conc) & + & \nabla \cdot (\dens \conc \dvel - \dens \disp \nabla \conc) = 0,
\end{eqnarray}
where the tensor field $\disp$ represents the molecular diffusion and the mechanical dispersion of the salt. For the velocity $\dvel$, Darcy's law is used:
\begin{eqnarray} \label {e_Darcy_vel}
 \dvel = - \frac{\perm}{\visc} (\nabla \pres - \dens \grav).
\end{eqnarray}
Here, $\pres = \pres (t, \mathbf{x}): [0, +\infty) \times \Do \to \mathbb{R}$ is the hydrostatic pressure and $\grav = (0, \dots, - g)^T \in \mathbb{R}^d$ the gravity vector.

In this work, we consider a special case of this general model. We assume that the porous medium is isotropic and characterized by the scalar permeability
\begin{eqnarray} \label {e_scalar_perm}
 \perm = K \mathbf{I},\;\; \text{where}\;\; K = K (\mathbf{x}) \in \mathbb{R},\;\;\text{and}\;\;\mathbf{I}\in\mathbb{R}^{d\times d}\;\text{is the identity matrix}.
\end{eqnarray}
For the density, we set
\begin {eqnarray} \label {e_lin_density}
 \dens (\conc) = \dens_0 + (\dens_1 - \dens_0) \conc,
\end {eqnarray}
where $\dens_0$ and $\dens_1$ are the densities of pure water and the brine, respectively. Note that $\conc \in [0, 1]$, where $\conc = 0$ corresponds to pure water and $\conc = 1$ to the saturated solution.  The viscosity is assumed to be constant. Finally, the dispersion is neglected:
\begin {eqnarray} \label {e_mol_diff}
 \disp = \poro D_m \mathbf{I},
\end {eqnarray}
where $D_m$ is the molecular diffusion coefficient.

Fields $\phi(\mathbf{x})$ and $K(\mathbf{x})$ will depend on the stochastic variables. Values of the other parameters used in this work are presented in Table \ref{tab:ElderParam}.
\begin{table}[b]
\begin{center}
 \caption{Parameters of the density-driven flow problem}
 \label{tab:ElderParam}
 \begin{tabular}{|l|l|l|} \hline
  Parameters & Values and Units & Description \\ \hline
  $\EXP{\phi}$ & 0.1 [-] & mean value of the porosity \\ \hline
  $D_m$ & $0.565\cdot 10^{-6}$ [$\mathrm{m}^2 \cdot \mathrm{s}^{-1}$] & molecular diffusion \\ \hline
  $\EXP{\perm}$ & $4.845\cdot 10^{-13}$ [$\mathrm{m}^2$] & mean value of the permeability \\ \hline
  $\grav$ & $9.81\be_z$ [$\mathrm{m} \cdot \mathrm{s}^{-2}$] & gravity \\ \hline
  $\rho_0$ & $1000$ [$\mathrm{kg} \cdot \mathrm{m}^{-3}$] & density of pure water \\ \hline
  $\rho_1$ & $1200$ [$\mathrm{kg} \cdot \mathrm{m}^{-3}$] & density of brine \\ \hline
  $\mu$ & $10^{-3}$ [$\mathrm{kg} \cdot \mathrm{m}^{-1} \cdot s^{-1}$] & viscosity \\ \hline
 \end{tabular}
\end{center}
\end{table}
%

%The problem (\ref {e_cont_eq}--\ref {e_mol_diff}) must be closed by the %specification for the boundary conditions for equations (\ref {e_cont_eq}) %and (\ref{e_tran_eq}), i.e. for $\conc$ and $\pres$, as well as the initial %conditions for $\conc$. 
In our numerical experiments, variants of the Elder problem \cite{Voss_Souza,Elder_Overview17} are considered. In it, concentrated solution intrudes the upper boundary into an aquifer filled with pure water. Due to the difference in the density, a complicated, unstable flow arises. This leads to a specific distribution of the mass fraction typically described as ``fingering'' \cite{Johannsen2003}. The time evolution of the mass fraction and the pressure in (\ref{e_cont_eq}--\ref{e_Darcy_vel}) is determined by the initial conditions for $\conc$. Note that the system (\ref {e_cont_eq}--\ref {e_Darcy_vel}) is nonlinear, and the model may have several stationary states \cite{Johannsen2003}. Any changes of the porosity and permeability fields may have an essential effect on the solution: The same initial and boundary data may correspond to principally different asymptotic behaviors. The influence of this factor can be estimated by the application of the stochastic models.

\begin{figure}[htbp!]
   \begin{subfigure}[b]{0.48\textwidth}
    \centering
     \caption{}
    \includegraphics[width=0.69\textwidth]{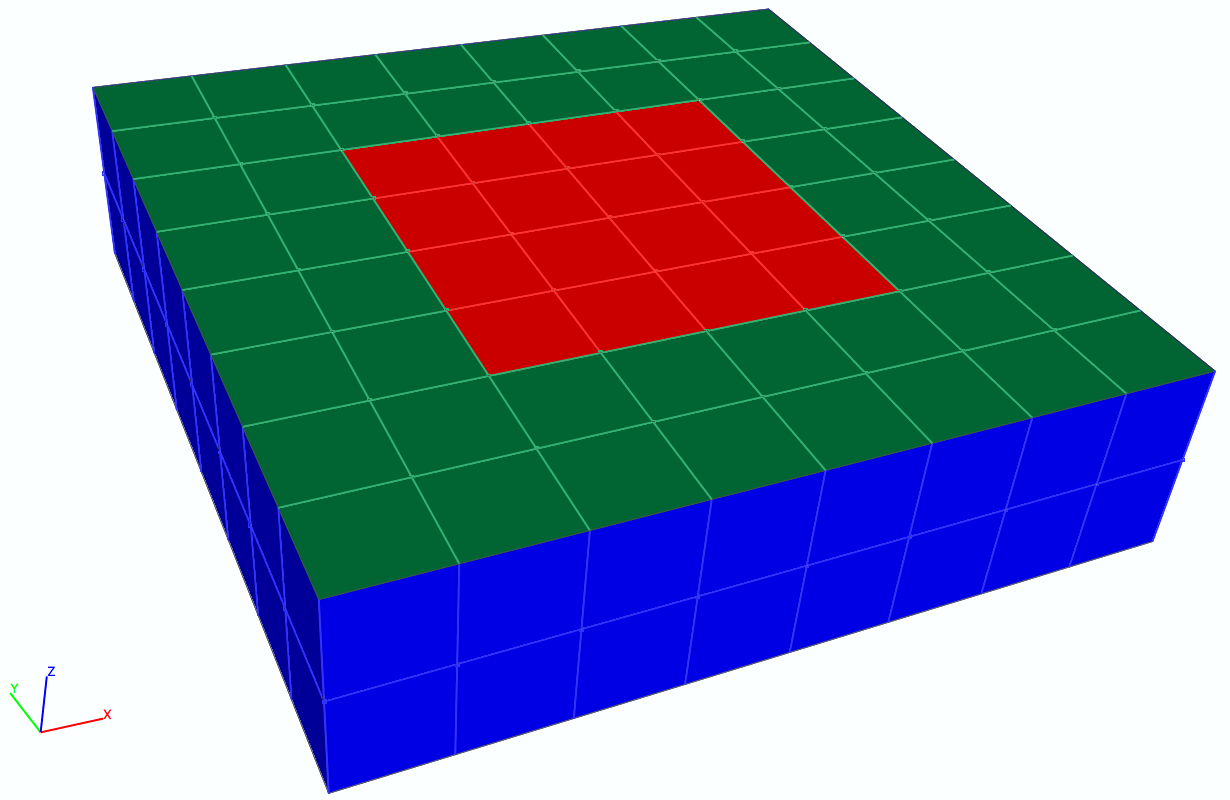}
     \label{fig:Elder3d-geom}
    \end{subfigure}
   \begin{subfigure}[b]{0.48\textwidth}
    \centering
     \caption{}
    \includegraphics[width=0.79\textwidth]{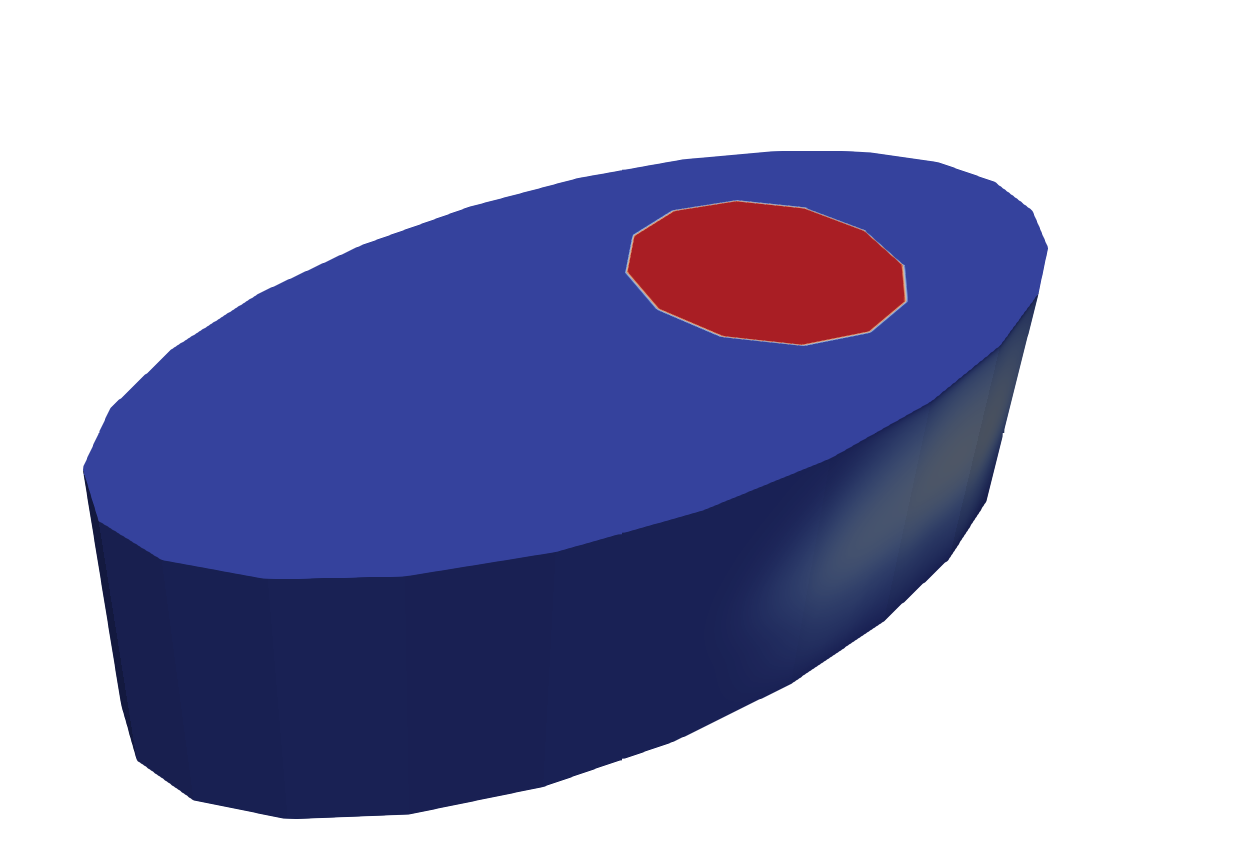} % Elder_ns_vtk/
      \label{subfig:Elder3d-cyl}
    \end{subfigure}
    \caption{Geometries of the domains used in the 3d simulations. Boundary conditions for the salt: $\conc = 1$ in the red area on the top and $\conc = 0$ otherwise.}
    \label{fig:Elder3d-scheme}
\end{figure}

We consider two spatial 3D domains, cf. Fig.~\ref{fig:Elder3d-geom}: A rectangular parallelepiped $\D = [0,600] \times [0,600] \times [0,150]$ $\mathrm{m}^3$ (Fig.~\ref{fig:Elder3d-geom}(a)) and an elliptic cylinder (Fig.~\ref{fig:Elder3d-geom}(b)) with $x\in [-300,300]$ (the  long side), $y\in [-150,150]$ (the short side) and $z\in [-150,0]$ (the depth).
%The 3D computational domains considered in this work are shown in Fig.~ \ref{fig:Elder3d-geom}. The domain on the left is a parallelepiped $\D = [0,600] \times [0,600] \times [0,150]$ $\mathrm{m}^3$. The number of cells in the presented experiments varied from $524{,}288$ till $4{,}194{,}304$. 
%The second domain is a cylinder, where $x\in [-300,300]$ (the  long side), $y\in [-150,150]$ (the short side) and $z\in [-150,0]$ (the depth). The polluted area is a circle with the center in $(-150,0,0)$, and the radius $100$. 
%The number of cells in the presented experiments varied from $241{,}152$ till $15{,}433{,}728$.
We impose the Dirichlet boundary condition $\sol = 1$ on the red area on the top of $\mathcal{D}$ ($z = 0$), and $\sol = 0$ on the green area on the top and on the other boundaries. For equation \refeq{e_cont_eq}, we impose the no-flux boundary conditions on the bottom and the vertical walls. The boundary conditions for $p$ differ for both geometries: For the parallelepiped, we impose the no-flux boundary conditions on all the boundaries. However, to fix $p$, we set $p = 0$ on the edges between the green and the blue faces. For the cylinder, we set $p = 0$ on the whole top boundary.
%!!!

For the spatial discretization of \refeq{e_cont_eq}--\refeq{e_Darcy_vel}, a vertex-centered FV scheme is used. Degrees of freedom (DoFs) for $\conc$ and $p$ are located in the grid nodes of the primary mesh. The domain in Fig.~\ \ref{fig:Elder3d-scheme}(a) is covered by a grid of $524{,}288$ hexahedra, so that the total number of the DoFs is $1{,}098{,}306$. For the domain in Fig.~\ref{fig:Elder3d-scheme}(b), we use an unstructured mesh of $241{,}152$ tetrahedra with $89{,}586$ DoFs.

Resolution of the spatial and temporal discretizations play an important role in the accuracy of the simulations. The distribution of the hydrogeological parameters of the porous medium should be represented correctly. Besides that, smoothing the flow on a coarse grid due to the numerical diffusion can lead to a loss of some phenomena and therefore reduce the accuracy of the uncertainty quantification. For this, the grid convergence of the simulations w.r.t.\ the spatial mesh has been tested. We ran our simulation on the twice regularly refined spatial mesh for geometry from Fig.~\ref{fig:Elder3d-scheme}(b). It consists of $15{,}433{,}728$ cells and has $2{,}822{,}552$ DoFs. The obtained solution was similar to the solution, computed on a coarser mesh mentioned above. Thus the further refinement does not change the behaviour of the solution essentially. To achieve reasonable computation times, parallelization of the evaluation of every realization is performed.
\subsection{Stochastic modeling of porosity, permeability and mass fraction}
\label{litv:subsec:PorosityVar}
Due to lack of knowledge or inaccurate measurements, the 
two input parameters porosity and permeability are unknown. 
One of the approaches to deal with this situation is to say that
both input parameters are uncertain. For example, we can introduce a prior distribution for these two values. We also can, based on expert knowledge, assume that the mean and the variance values are known.

In the following we model the porosity ($\poro(\bx)$) by a random field ($\poro(\bx,\btheta)$). We also assume that the permeability ($\perm$) is a function of $\poro$ as shown in \refeq{e_perm_of_poro} and \refeq{e_perm_Kozeny_Carman}:
\begin {eqnarray} \label {e_perm_of_poro}
 \perm = K(\poro) \mathbf{I}.
\end {eqnarray}
Since parameters are uncertain the solution $\mysol$ is uncertain too, and this solution
is a function of $\poro$ and $\perm$.
Here we assume $\perm$ to be isotropic.
Usually the distribution of $\poro(\bx,\thetab)$, $\bx\in \Do$, is non-trivial or even very complicated. Therefore, it is comfortable from the numerical point of view, to represent this complicated random field in some easy to compute basis, e.g. with Gaussian random variables. Often, many Gaussian random variables are required.
We introduce $\btheta=(\theta_1,...,\theta_M,...)$, where $\theta_i$ are Gaussian random variables.

The dependence (\ref{e_perm_of_poro}) is specific for every material and there is no a general law. In our experiments, we use a Kozeny-Carman-like dependence \cite{Costa_2006}:
\begin{eqnarray} 
\label{e_perm_Kozeny_Carman}
 K (\poro) = \kappa_{KC} \cdot \dfrac {\poro^3} {1 - \poro^2},
\end {eqnarray}
where the scaling factor $\kappa_{KC}$ is chosen to satisfy the equality $K(\EXP{\poro}) \mathbf{I} = \EXP{\mathbf{K}}$. Further parameters of the standard Elder problem are listed in Table \ref{tab:ElderParam}.\\
%
%\subsection{Quantities of interest}
%
\subsection{Solution of the stochastic flow model}
\label{sec:NumFlow}
We start with sampling random variables $\{\btheta_i\}$ and computing the realisations of porosity $\poro_i(\bx)=\poro(\bx,\btheta_i)$ and permeability $\perm_i(\bx)=\perm(\bx,\btheta_i)$, and then solving equations (\ref{e_cont_eq}--\ref{e_tran_eq}). To compute the solution we used plugin d${}^3$f of the simulation framework \myug,
\cite{ug4_ref1_2013,ug4_ref2_2013}. This framework presents a flexible tool for the numerical solution of non-stationary and nonlinear systems of PDEs, and is parallelized and highly scalable.
Equations~\ref{e_cont_eq}--\ref{e_tran_eq} are discretized by a vertex-centered finite-volume method on
unstructured grids in the geometric space, as presented in \cite{Frolkovic-MaxPrinciple}. In particular, the consistent
velocity approach is used for the approximation of the Darcy velocity (\ref{e_Darcy_vel}),
\cite{Frolkovic-ConsVel,Frolkovic-Knaber-ConsVel}. The implicit Euler scheme is used for the time discretization.
The implicit time-stepping scheme is chosen for its unconditional stability, as the velocity depends on the
unknowns of the system. This is especially important for variable coefficients as it is difficult to predict the
range of the variation of the velocity in the realizations.

In this discretization, we obtain a large sparse system of nonlinear algebraic equations in every time step, and we solve it using the Newton method with a line search method. The sparse linear system, which appears in the nonlinear
iterations, is solved by the BiCGStab method \cite {Templates} with the multigrid preconditioning (V-cycle,
\cite{Hackbusch85}), which proved to be very efficient in this case. In the multigrid cycle, we use the
ILU${}_\beta$-smoothers \cite{Hackbusch_Iter_Sol}.
\section{Stochastic Modeling and Methods}
\label{litv:sec:StochModel}
%
%\subsection{\gPC-based surrogate}
%\label{litv:sec:methods}

%This surrogate representation is built step by step
%or sample by sample, and can already be used for each new sample. In case we are
%sampling a solution of an SPDE, this allows us to reduce the number of necessary
%samples namely in case the solution is already well-represented by the low-rank
%tensor approximation. This can be easily checked by evaluating the residuum of
%the PDE with the approximate solution.

Let $\btheta=(\theta_1,\ldots,\theta_M)$ be a vector of independent random variables.
Let ${S}:=(\vTheta, \mathcal{B}, \mathbb{P})$ be a probability space.
Here $\vTheta$ denotes a sample space, $\mathcal{B}$ a $\sigma$-algebra on $\vTheta$, and $\mathbb{P}$ a probability measure on $(\vTheta, \mathcal{B})$. 
Random variable $\theta:\vTheta \rightarrow \mathbb{R}^M$ with finite variance can be represented in terms of the polynomial chaos expansion, e.g. as a multi-variate Hermite polynomial of Gaussian RVs, \cite{wiener38}.
Alternatively to Hermite polynomials, other orthogonal polynomials (Legendre, Chebyshev, Laguerre) \cite{Askey85, xiuKarniadakis02a, Xiu}, or splines, for instance, can be used.
We assume $\{\Pol_{\beta} \}_{\beta \in \mathcal{J}}$ is a basis of $S=L_2(\vTheta, \mathbb{P})$. The cardinality of the multi-index set $\mathcal{J}$ is infinite, therefore, we truncate it and obtain a finite set $\mathcal{J}_{M,p}$ and a subspace $S_{M,p}$. The solution $\mysol(t,\bx,\btheta)$ lies in the tensor product space $L_2(D)\otimes L_2(\vTheta)\otimes L_2([0,T])$, where $[0,T]$ is a time interval, where the initial problem is solved. After a discretization, we obtain $\bsol \in \mathbb{R}^{n}\otimes {S}_{M,P}\otimes \mathbb{R}^{n_t}$, where $n$ is the number of basis functions in the physical space, $\vert \mathcal{J}_{M,p} \vert$ dimension of the stochastic space, and $n_t$ dimension of the temporal space.

The empirical mean value and the variance of $\mysol(t,\bx,\btheta)$ can be computed as follow
%an integral over the multidimensional domain $\vTheta$:
\begin{equation} \label{eq:pdestat}
  \overline{\mysol}(t,\bx)=\EXP{\mysol(t,\bx,\btheta)} = \int_\vTheta \mysol(t,\bx,\btheta)\,\mathbb{P}(\di \btheta)\approx \sum_{i=1}^{N_q} w_i \mysol(t,\bx,\btheta_i),
\end{equation}
with some quadrature points $\btheta_i$ and weights $w_i$.
Here $\mysol(t,\bx,\btheta_i)$ is the solution evaluated at fixed point
$\btheta=\btheta_i \in \mathbb{R}^M$.
The sampled variance is computed as follow
\begin{equation} \label{eq:numer_var}
  \var{\mysol}(t,\bx)\approx \sum_{i=1}^{N_q} w_i \left ( \overline{\mysol}(t,\bx) - \mysol(t,\bx,\btheta_i) \right)^2.
\end{equation}
Other well-known methods to compute the mean and the variance are quasi-Monte Carlo, Monte Carlo and multi-level Monte Carlo methods.
\subsection{Generalized Polynomial Chaos Expansion}
\label{litv:sec:gPCE}
As outlined in \cite{Dongbin, xiuKarniadakis02a, wiener38}, the random variable (field) $\Psii(t,\bx,\thetab)$ can be represented as a series of multivariate Legendre polynomials in uncorrelated and independent uniformly distributed random variables
 $\thetab=(\theta_1,...,\theta_j,...)$:
\begin{equation}
\label{eq:PCEdecom}
\Psii(t,\bx,\thetab)=\sum_{\beta \in \mathcal{J}}\Psii_{\beta}(t,\bx)\Pol_{\beta}(\thetab),
\end{equation}
where $\{\Pol_{\beta}\}$ is a multivariate Legendre basis, $\beta = (\beta_1 ,..., \beta_j ,...)$ is a multiindex and $\mathcal{J}$ is a multiindex set. 

This expansion is called the \textit{generalized polynomial chaos expansion} (\gPC) (see more in the Appendix or in  
\cite{wiener38, xiuKarniadakis02a, Xiu, ernst_mugler_starkloff_ullmann_2012}). 
For numerical purposes, we truncate this infinite \gPC series and obtain a finite series approximation, i.e., we keep only $M$ random variables, $\thetab=(\theta_1,...,\theta_M)$, and limit the maximal order of the multi-variate polynomial $\Pol_{\beta}(\thetab)$ by $p$. We denote the new multiindex subset as $\mathcal{J}_{M,p} \subset \mathcal{J}$. Then the unknown random field $\Psii(t,\bx,\thetab)$ can be approximated as follows:
\begin{equation}
\label{eq:PCEdecom2}
\Psii(t,\bx,\thetab)\approx \widehat{\Psii}(t,\bx,\vTheta)=\sum_{\beta \in \mathcal{J}_{M,p}}\Psii_{\beta}(t,\bx)\Pol_{\beta}(\thetab).
\end{equation}
Here $\Pol_{\beta}(\btheta)$ are multi-variate Legendre polynomials defined as follows
\begin{equation*}
  \label{eq:B:def-mult-Legendre}
  \Pol_\beta(\thetab) := \prod_{j=1}^{M} \pol_{\beta_j}(\theta_j);
  \quad \forall \thetab\in\mathbb{R}^M,\, \beta \in\mathcal{J}_{M,p},\;\sum_{j=1}^M\beta_j \leq p, 
\end{equation*}
%$\beta = (\beta_1,\ldots,\beta_j,\ldots)
%\in \mathcal{J}$
$\pol_{\beta_j}(\cdot)$ are Legendre monomials, defined in \refeq{eq:LegendreMonom}. The scalar product is 
\begin{equation}
\label{eq:spro}
\langle \Pol_{\alpha},\Pol_{\beta} \rangle_{L_2(\vTheta)} = \EXP{\Pol_{\alpha}(\thetab)\Pol_{\beta}(\thetab)}=\int_{\Theta}\Pol_{\alpha}(\thetab)\Pol_{\beta}(\thetab)d\mathbb{P}{\thetab}=Q_{\alpha}\delta_{\alpha\beta},
\end{equation}
where $Q_{\alpha}=\EXP{\Pol_{\alpha}^2}=q_{\alpha_1}\cdot...\cdot q_{\alpha_d}$, $q_{\alpha_j}=\langle \pol_{\alpha_j},\pol_{\alpha_j} \rangle=1/(2\alpha_j +1)$, are the normalization constants and $\delta_{\alpha\beta}=\delta_{\alpha_1\beta_1}\cdot ... \cdot \delta_{\alpha_M \beta_M}$ is the $M$-dimensional Kronecker delta function.

Decomposition \eqref{eq:PCEdecom} can be understood as a response surface for $\Psii(t,\bx,\thetab)$. As soon as the response surface is built, the value $\Psii(t,\bx,\btheta)$ can be evaluated for any $\btheta$ almost for free (only by evaluating the polynomial (\ref{eq:PCEdecom})). It can be very practical if $10^6$ samples of $\Psii(t,\bx,\btheta)$ are needed, for example.
Since Legendre polynomials are $L_2$-orthogonal, the coefficients $\Psii_{\beta}(t,\bx)$ can be computed by projection: 
\begin{equation}
\label{eq:PCEcoef}
\Psii_{\beta}(t,\bx)=\frac{1}{\langle \Pol_{\beta},\Pol_{\beta} \rangle}\int_{\Theta}\widehat{\Psii}(t,\bx,\btheta)\Pol_{\beta}(\btheta)\,\mathbb{P}(\di \btheta)\approx \frac{1}{\langle \Pol_{\beta},\Pol_{\beta} \rangle}\sum_{i=1}^{N_q}\Pol_{\beta}(t,\thetab_i)\widehat{\Psii}(t,\bx,\thetab_i)w_i,
\end{equation}
where $w_i$ are weights and $\thetab_i$ quadrature points, defined, for instance, by a Gauss-Legendre integration rule.
Once all \gPC coefficients are computed we substitute them into \refeq{eq:PCEdecom}. For estimates of truncation and aliasing errors, see \cite{ConradMarzouk13}.
There are different strategies to truncate the multi-index set $\mathcal{J}$ to $\mathcal{J}_{M,p}$.
In the multi-index set $\mathcal{J}_{M,p}$ we keep only $M$ RVs, i.e., $\thetab=(\theta_1,\ldots,\theta_M)$ and $p_i$ is the monomial degree w.r.t. the variable $\theta_i$.
For the total multi-variate polynomial degree $p$ we pose one of the following condition
a) $\prod_{i=1}^M p_i \leq p$ or b) $\sum_{i=1}^M p_i \leq p$ or c) $p_i\leq p$, $\forall\; i=1..M$.

\begin{remark}
If only the variance and the mean values are of interest, we recommend to use the collocation method with some sparse or full-tensor grid. If some additional statistics are required, such as sensitivity analysis, parameter identification, or computing cdf or pdf, we recommend to compute the gPC expansion.
\end{remark}

Quadrature grids for computing multi-dimensional integrals can be constructed from products of 1D integration rules. To keep computational costs to a minimum, a sparse grid technique can be used for certain class of smooth functions \cite{Griebel_Bungartz}. To compute \refeq{eq:PCEcoef}, we apply the well-known Clenshaw Curtis quadrature rules \cite{CC60}. 
%Additionally to the full tensor grid, we use Smolyak's sparse grid procedure to construct a cheap sparse grid approximation. 
%

%
%
%
%\subsection{Legendre polynomials}
Since $\Pol_{\bzero}(\btheta)=1$, it becomes apparent that the mean value is the first \gPC  coefficient 
\begin{equation}
\label{eq:PCEcoef2}
\overline{\sol}:=\EXP{\Psii(t,\bx,\btheta)}=\Psii_{\bzero}(t,\bx)=
\frac{1}{\langle \Pol_{\bzero},\Pol_{\bzero} \rangle}\int_{\Theta}{\Psii}(t,\bx,\btheta)\Pol_{\bzero}(\btheta)\,\mathbb{P}(\di \btheta)
\approx \sum_{i=1}^{N_q}{\Psii}(t,\bx,\thetab_i)w_i, 
\end{equation}
$\bzero=(0,\ldots,0)$.
The variance is the sum of squared \gPC coefficients \cite{xiuKarniadakis02a, GhanemBook17,Knio10}:
\begin{align}
\label{eq:var}
\var{\Psii(t,\bx,\btheta)}&=\EXP{\left(\Psii(t,\bx,\btheta)-\overline{\Psii}\right) \otimes \left(\Psii(t,\bx,\btheta)-\overline{\Psii}\right)}\\
%&=\sum_{\beta \in \mathcal{J}, \beta>0}\sum_{\gamma \in \mathcal{J}, \gamma>0} %\EXP{\Pol_{\beta}(\thetab)\Pol_{\gamma}(\thetab)}\Psii_{\beta}(t,\bx)\otimes %\Psii_{\gamma}(t,\bx)\\
&=\sum_{\gamma>0}\Vert \Pol_{\gamma} \Vert^2 \Psii_{\gamma}(t,\bx)\otimes \Psii_{\gamma}(t,\bx)=
\sum_{\gamma>0}\prod_{k=1}^M\frac{1}{2\gamma_k+1} \Psii_{\gamma}(t,\bx)\otimes \Psii_{\gamma}(t,\bx),
\end{align}
where $\otimes$ is the Kronecker product.
\subsection{Computing probability density functions}
After we calculated the \gPC surrogate $\widehat{\sol}(t,\bx,\thetab)$ in \refeq{eq:PCEdecom2}, we can compute the probability density functions, exceedance probabilities, and quantilies in the selected points.
For this, one should evaluate the multi-variate polynomial $\widehat{\sol}(t,\bx,\thetab)$ in a sufficiently large number of points $\thetab$
%Approximation of the probability density function in a point $(t^{*},\bx^{*})$ is computed by sampling the multivariate polynomial on the right-hand side in \refeq{eq:PCEpoint}. 
\begin{equation}
\label{eq:PCEpoint}
\widehat{\mysol}(\thetab):=\widehat{\mysol}(t^{*},\bx^{*},\thetab)=\sum_{\beta \in \mathcal{J}_{M,p}}\mysol_{\beta}(t^{*},\bx^{*})\Pol_{\beta}(\thetab).
\end{equation}
%These evaluations are of low coputing cost.
The random variable $\widehat{\mysol}(\thetab)$ can be sampled, e.g., $N_s=10^6$ times (no additional extensive simulations are required). The obtained sample can be used to evaluate the required statistics, and the probability density function. The exceedance probability for some threshold $c^{*}$ can be estimated as follows:
\begin{equation}
    P(\mysol>\mysol^{*})\approx \frac{\#\{\widehat{\mysol}(\thetab_i):\; \widehat{\mysol}(\thetab_i)>\mysol^{*}, \; i=1,\ldots,N_s\}}{N_s}.
\end{equation}
Having a sufficiently large sample set makes it possible to estimate quantiles \cite{LitvElder2D}.
%We remind that for a random variable $\xi$ and for a fixed value $\alpha\in (0,1)$, the $\alpha$-quantile is the number $\xi_{\alpha}\in \mathbb{R}$ such that $\mathbb{P}(X\leq \xi_{\alpha})\geq \alpha$.
%
%
\subsection{Numerical errors}
\label{sec:trunc}
Applying \gPC approximation, we introduce the truncation and approximation errors \cite{CONSTANTINE12, Sinsbeck_2015,ConradMarzouk13, ernst_mugler_starkloff_ullmann_2012}.
Truncating \gPC coefficients $\{\Psii_{\beta}:\;\beta \in \mathcal{J}_{c}\}$, we introduce the truncation error
\begin{equation}
\label{tr_err}
\err_t=\Vert \Psii(t,\bx,\thetab)- \widehat{\Psii}(t,\bx,\thetab)\Vert =
\Big \| \sum_{\beta \in \mathcal{J}_{c}}\Psii_{\beta}(t,\bx)\Pol_{\beta}(\thetab) \Big \|,  \quad \mathcal{J}_{M,p} \cup \mathcal{J}_{c}=\mathcal{J}.  
\end{equation}
Additionally, we approximate the coefficients $\Psii_{\beta}(t,\bx)$ by $\widehat{\Psii}_{\beta}(t,\bx)$. This introduces the approximation error
\begin{equation}
\label{eq:appr_err}
\err_a=\Big \Vert
\sum_{\beta \in \mathcal{J}_{M,p}}\Psii_{\beta}(t,\bx)\Pol_{\beta}(\thetab)
-
\sum_{\beta \in \mathcal{J}_{M,p}}\widehat{\Psii}_{\beta}(t,\bx)\Pol_{\beta}(\thetab) \Big \Vert
= \Big \Vert \sum_{\beta \in \mathcal{J}_{M,p}}(\widehat{\Psii}_{\beta}(t,\bx)-\Psii_{\beta}(t,\bx))\Pol_{\beta}(\thetab) \Big \Vert
\end{equation}
The sum of both errors is
\begin{equation}
\label{eq:sum_errs}
\err_t+\err_a=\underbrace{\Big \Vert \sum_{\beta \in \mathcal{J}_c}{\Psii}_{\beta}(t,\bx)\Pol_{\beta}(\thetab) \Big \Vert}_\text{truncation error} +\underbrace{\Big \Vert \sum_{\beta \in \mathcal{J}_{M,p}}(\Psii_{\beta}(t,\bx)-\widehat{\Psii}_{\beta}(t,\bx))\Pol_{\beta}(\thetab) \Big \Vert}_\text{approximation error}.
\end{equation}
\begin{remark}
Spatial resolution plays an important role in the simulations. Thus, the distribution of the parameters of the porous medium should be represented correctly. Furthermore, smoothing the flow on a coarse grid due to numerical diffusion can lead to a loss of some phenomena and therefore reduce the accuracy of the uncertainty quantification. Thus we must cover $\Do$ with sufficiently fine grids and use small time steps. To achieve reasonable computation times, this also assumes parallelization of the evaluation of every realization.
\end{remark}
%\newpage

%
\section{Parallelization}
\label{sec:parallelization}

The computation of the mean value, variance, and other statistics are done in parallel.
There are two levels of parallelization. First, we compute the solution in all quadrature or sampling points in parallel no communication is needed). Second, the solution in each single quadrature point is computed also in parallel (communication is needed).

We performed computations on a Cray XC40 parallel cluster Shaheen II provided by the King Abdullah University of Science and Technology (KAUST) in Saudi Arabia. Shaheen II has 6174 nodes with 2 Intel Haswell 2.3 GHz CPUs per node, and 16 processor cores per CPU. Every node has 128 GB RAM. 
%The nodes communicate via the Cray Aries interconnect network with Dragonfly topology. The system is equipped with the %Sonexion Lustre 2000 storage appliance.

On the first level, we allocated, depending on the experiment, up to 600 computing nodes. It allowed us to compute (up to) 600 \QMC simulations simultaneously. All these simulations are independent, and, therefore, no communication between nodes was required.
The total computing time is equal to the time, required for computing one simulation. On the second level, on each node, we allocated 32 computing cores to resolve each single simulation. Thus, in total, we run our code on up to $600 \times 32=19200$ parallel cores. If a simulation is large or requires a lot of memory, we can allocate more than one node for it. \myug allows allocating all available nodes (6174) for solving just one expensive simulation. Theoretically, \myug allows using all available nodes and cores, i.e., $6174\times 32=197{,}568$ cores.

Parallelization of the framework \myug is based on the distribution of the spatial grid between the processes and is implemented using MPI \cite{ug4_ref1_2013}. As the computation of the realizations involves the nonlinear iterations and the solution of the large sparse linear systems (see Section \ref{sec:NumFlow}), quite intensive communication is required. 
Shaheen II used the SLURM system to manage parallel jobs, and the job arrays were used for the concurrent computations of the realizations. This was the most time-consuming phase, which required hours of computing time.
\section{Numerical Experiments}
\label{litv:sec:numerics}
\subsection{Solving the deterministic problem}
\label{sec:Num}
The deterministic solver \myug uses the multigrid solver with a hierarchy of $L$ nested spatial grids in $\D$. There are $N_{\ell}$ nodes on each level $\ell=1,\ldots,L$. Each grid is distributed over all available processors on the node (32 in the tests below).

Equations~\ref{e_cont_eq}--\ref{e_tran_eq} are discretized by a vertex-centered finite-volume method on unstructured grids in space, as presented in \cite{Frolkovic-MaxPrinciple}. In particular, we use the ``consistent velocity'' for the approximation of the Darcy velocity \eqref{e_Darcy_vel}, cf.\ \cite{Frolkovic-ConsVel,Frolkovic-Knaber-ConsVel}. The implicit Euler scheme is used for the time discretization. The choice of the implicit time stepping scheme is motivated by its unconditional stability, as the velocity depends on the unknowns of the system.

After discretization one obtains a large sparse system of nonlinear algebraic equations in every time step. It is solved by the Newton method with a line search method. The sparse linear system, which appears in the nonlinear iterations, is solved by the BiCGStab method (cf.\ \cite {Templates}) with the multigrid preconditioning (V-cycle, cf.\ \cite{Hackbusch85}) which proved to be very efficient in this case. In the multigrid cycle, we use the ILU${}_\beta$-smoothers \cite{Hackbusch_Iter_Sol}. The matrices are assembled as the Jacobians of the discretized nonlinear systems for each grid. On each step of the Newton method the inversion of the Jacobian is computed. Developers of \myug also efficiently solve the problem of optimal distribution of the hierarchy of grids onto parallel processes.

%%%%%%%%%%%%%%%%%%%%%%%%%%%%
%
%

\textbf{Software.}  
\myug is a novel flexible software system for simulating PDE based models on high performance parallel clusters \cite{ug4_ref1_2013,ug4_ref2_2013}. The software toolbox \myug has been developed since the early 1990's. This software framework has successfully been applied to a variety of problems such as density-driven and thermohaline flow in porous media, Navier-Stokes equation, drug diffusion through human skin, level-set methods.
The aim of the software framework \myug is to solve the above mentioned problems on adaptive hybrid grid hierarchies in 2D and 3D dimensions. Special care is taken to parallel performance issues, scalability, linear algebra algorithms, and data structures. In addition, \myug is specially designed for user-friendly handling and functionality is made available to non-programming users by the usage of scripts and graphical representations.
The numerical solution of the systems (\ref {e_cont_eq}--\ref {e_tran_eq}) has been performed using the D3F plugin of the simulation software package \myug, cf.\ \cite{ug4_ref1_2013,ug4_ref2_2013}.
\subsection{3D computational domains}
The geometry of real-life 3D computational domains can be very complicated. Discretisation of such domains (reservoirs) may require non-trivial 3D meshes. This could be very challenging in the stochastic context, which we have here. Therefore, we start with two relatively trivial reservoir geometries: 1) a parallelepiped $\mathcal{D}=[0,600]\times [0,600]\times[0,150]$ (see Fig.~\ref{fig:Elder3d-geom}); 2) an
elliptic cylinder with the semi-axes $600$ and $300$ for the base ellipse and height $150$. 
%where the longest side is $x\in [-300,300]$, the short side $y\in [-150,150]$ and the depth 
Thus, the longest radius is $300$, $x\in[-300,300]$, the shortest radius is $150$, $y\in[-150,150]$, the height $z\in [-150,0]$ and the center is in the point $(0,0,-75)$ (see Fig.~\ref {subfig:Elder3d-cyl}). 
The polluted area is a circle with the center in $(-150,0,0)$, and the radius $100$.
%
% \paragraph{1. Parallelepiped.}
% We assume that the reservoir has three layers. The porosity coefficient is defined as follow
% \begin{align*}
% \phi(\bx,\thetab(\omega))&=0.1+ \exp(\theta_1\sin(\pi x/600)+
% \theta_2\sin(\pi y/600)+ \theta_3\sin(\pi z/150)+
% \theta_1\sin(\pi x/600)+ \\
% &+\theta_1\sin(\pi x/600)\sin(\pi y/600)
% +\theta_2\sin(\pi x/600)\sin(\pi z/150)+\theta_3\sin(\pi y/600)\sin(\pi z/150)).
% \end{align*}
% Figure~\ref{fig:3D_poro} shows a realization of the porosity field.
% %We computed gPCE coefficients and then the mean and the variance of the mass fraction. 
% Five isosurfaces of \textcolor{red}{WHICH SOLUTION? PROBABLY WILL REMOVE IT}the solution are shown in Figs.~\ref{subfig:refsol-up} (top view) and \ref{subfig:refsol-bottom} (bottom view).
% For this experiment we used the full tensor Gauss-Legendre quadrature grid with 27 points. The finest spatial mesh contained $1{,}098{,}306$ nodes. 

\subsection{Parallelepiped with 3 RVs}
We consider an experiment with three random variables $\xib=(\xi_1,\xi_2,\xi_3)$ and with the following porosity coefficient
\begin{equation}
\label{eq:poro-paral3RVS}
\poro(\bx,\xib):=0.1+0.01\left (\xi_1\sin\frac{x\pi}{600} +\xi_2\sin\frac{y\pi}{600} +
\xi_3\sin\frac{z\pi}{150} +
\xi_1\sin\frac{x\pi}{600}\sin\frac{y\pi}{600}+
\xi_2\sin\frac{x\pi}{600}\sin\frac{z\pi}{150}
\right),
\end{equation}
where $\bx=(x,y,z)\in [0,600]\times[0,600]\times[0,150]$ (parallelepiped).

\begin{figure}[htbp!]
\center
    \begin{subfigure}[b]{0.47\textwidth}
    \centering
      \caption{}
  \includegraphics[width=0.99\textwidth]{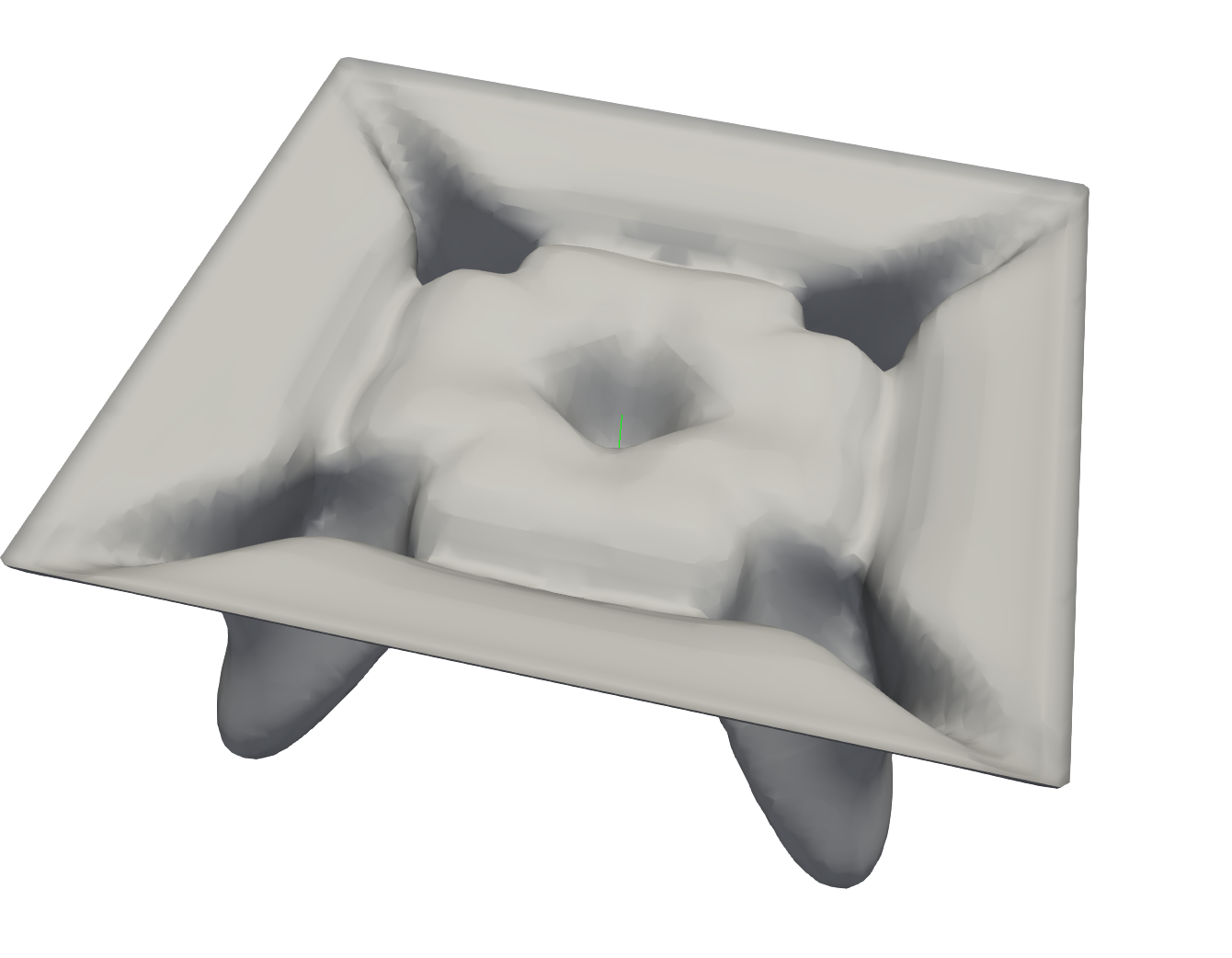}
     \label{subfig:contMC_mean}
    \end{subfigure}
  \begin{subfigure}[b]{0.47\textwidth}
    \centering
      \caption{}
  \includegraphics[width=0.99\textwidth]{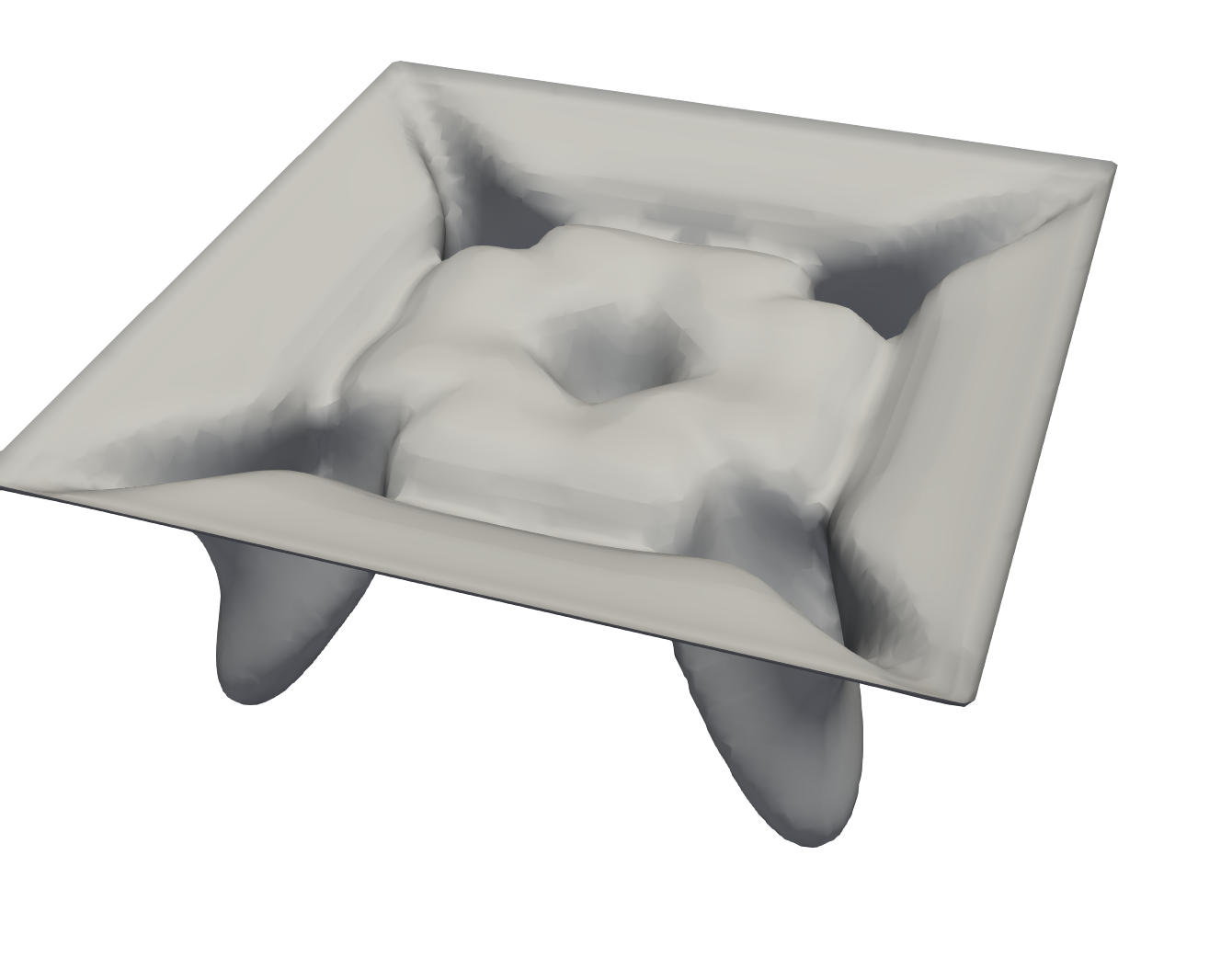}
     \label{subfig:contPC_mean}
    \end{subfigure}
\caption{Isosurface $\{\bx:\; \overline{\sol}(\bx)=0.5\}$ computed using (a) 200 \QMC simulations and (b) \gPC response surface of order 4. The upper square edge corresponds to the boundary of the red area in Fig.\ \ref{fig:Elder3d-scheme}(a)}
\label{fig:3D_ave_countur}
\end{figure}

In Fig.~\ref{fig:3D_ave_countur}, we visualize isosurfaces $\{\bx:\; \overline{\sol}(\bx)=0.5\}$ after 1500 time steps ($0.005\cdot 1500=7.5$ years), computed using (a) 200 \QMC simulations (Halton sequence) and (b) \gPC response surface of order 4. All \gPC coefficients were computed via a 3D Gauss-Legendre quadrature rule over the domain $[-1,1]^3$ with 125 quadrature points.

Note that in the present problem setting, the flow field is very complicated. The main process tends to represent the same behaviour as in the 2D version of the Elder problem. The heavy concentrated solution from the top boundary falls down and replaces the pure water located initially at the bottom. This pure water is pressed up. As the diffusion is very small, the transport of the salt by these opposite flows creates a non-trivial pattern of the concentration field typically described as ``fingering''. The particular pattern depends on the profile of the boundary conditions for $\conc$ and the initial flow. In our case, we see 5 ``fingers'' (4 at the corners and one at the center).
As we see, the isosurfaces $\{\bx:\; \overline{\sol}(\bx)=0.5\}$ are almost identical for both the methods.

\begin{figure}[htbp!]
\center
    \begin{subfigure}[b]{0.32\textwidth}
    \centering
      \caption{}
  \includegraphics[width=0.99\textwidth]{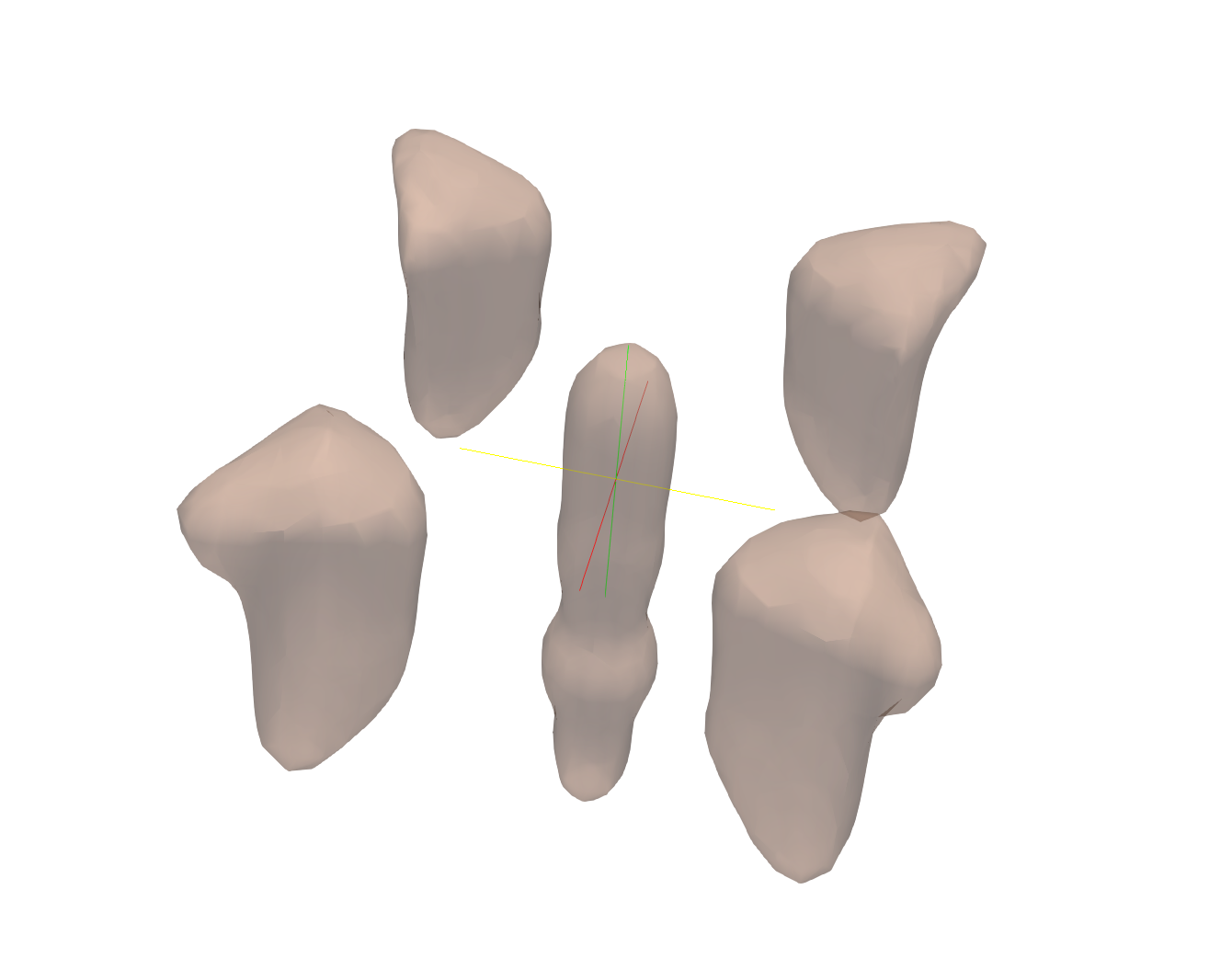}
     \label{subfig:Paral3D_MC200_ts1500_var_Countur}
    \end{subfigure}
  \begin{subfigure}[b]{0.32\textwidth}
    \centering
      \caption{}
  \includegraphics[width=0.99\textwidth]{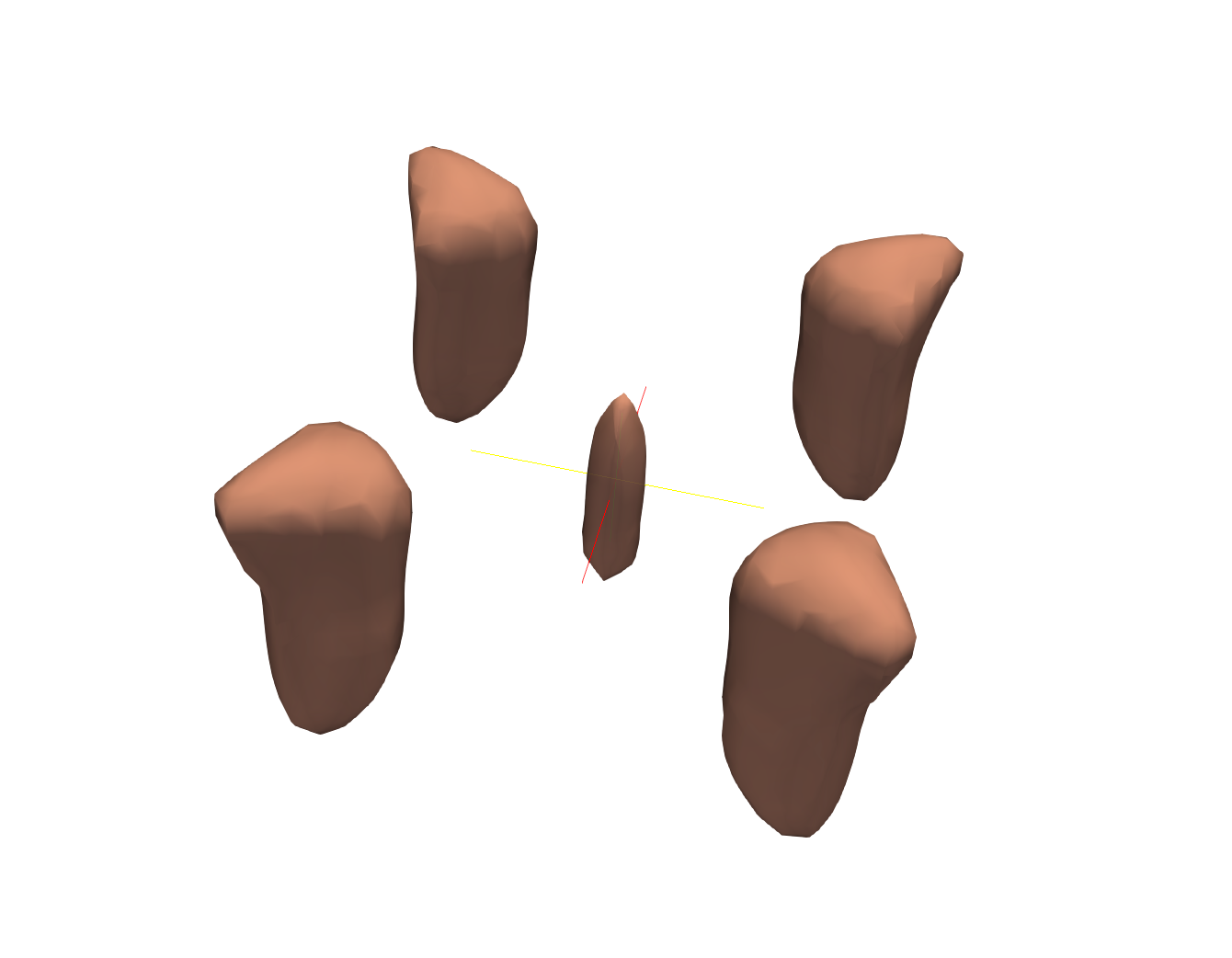}
     \label{subfig:Paral3D_GLF_PC35_ts1500_var_Countur}
    \end{subfigure}
  \begin{subfigure}[b]{0.32\textwidth}
    \centering
      \caption{}
  \includegraphics[width=0.99\textwidth]{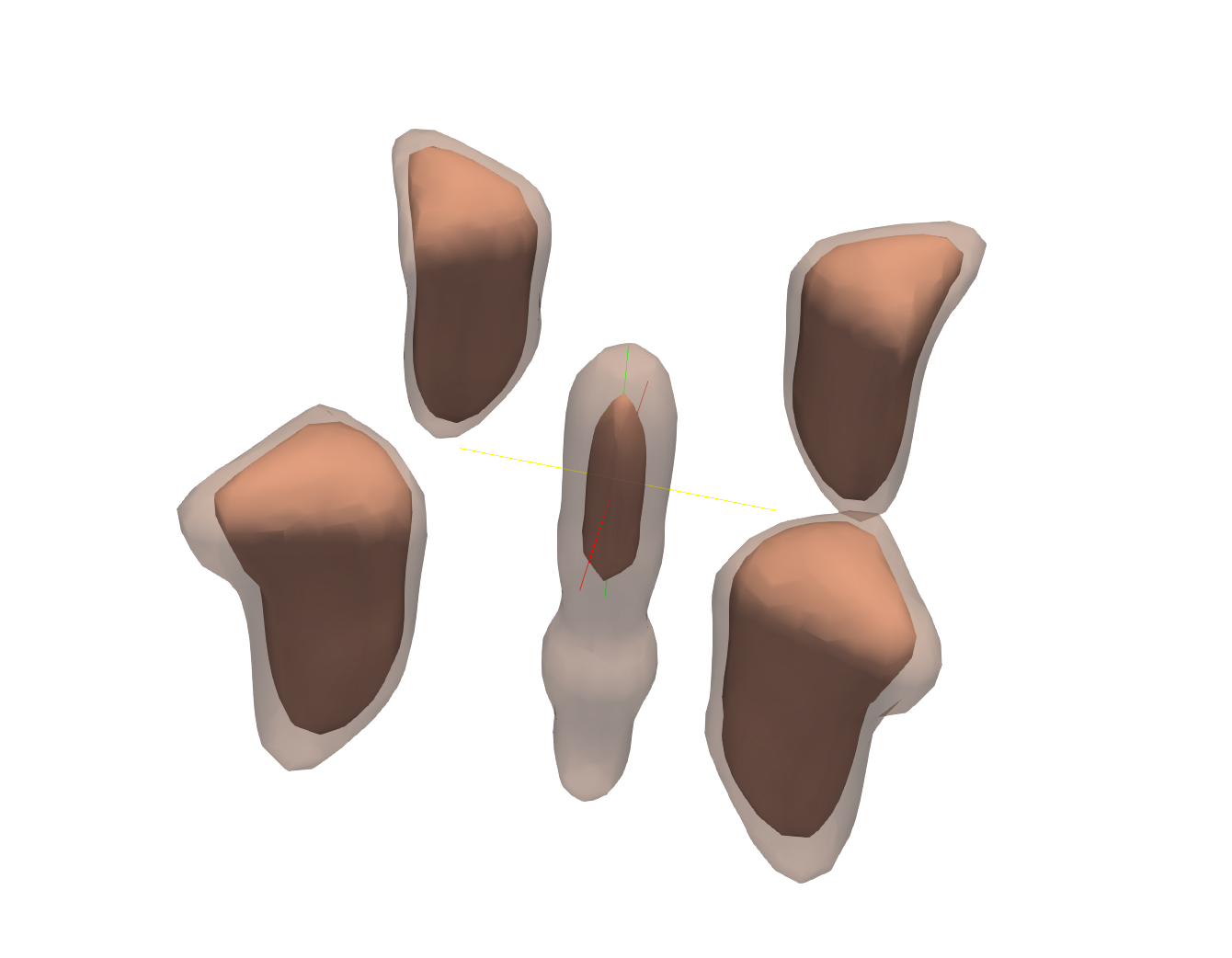}
     \label{subfig:cont_var_both}
    \end{subfigure}
\caption{Isosurface $\{\bx:\;\var{\sol}=0.05\}$, computed via a) \QMC and b) \gPC response surface of order 4, and c) comparison of both isosurfaces.}
\label{fig:3Dparal_var_countur}
\end{figure}

The same holds for the variance. Fig.~\ref{fig:3Dparal_var_countur} presents the isosurface $\{\bx:\;\var{\sol}(\bx)=0.05\}$ after these 1500 time steps for (a) 200 \QMC simulations (Halton sequence), (b) \gPC response surface of order 4 and (c) comparison of the isosurfaces from (a) and (b). We see that positions, numbers, and shapes of all fingers are similar. The isosurfaces, computed via \QMC, are slightly larger, especially one in the middle (see Fig.~\ref{subfig:cont_var_both}).

\begin{figure}[htbp!]
\center
    \begin{subfigure}[b]{0.48\textwidth}
    \centering
      \caption{}
  \includegraphics[width=0.99\textwidth]{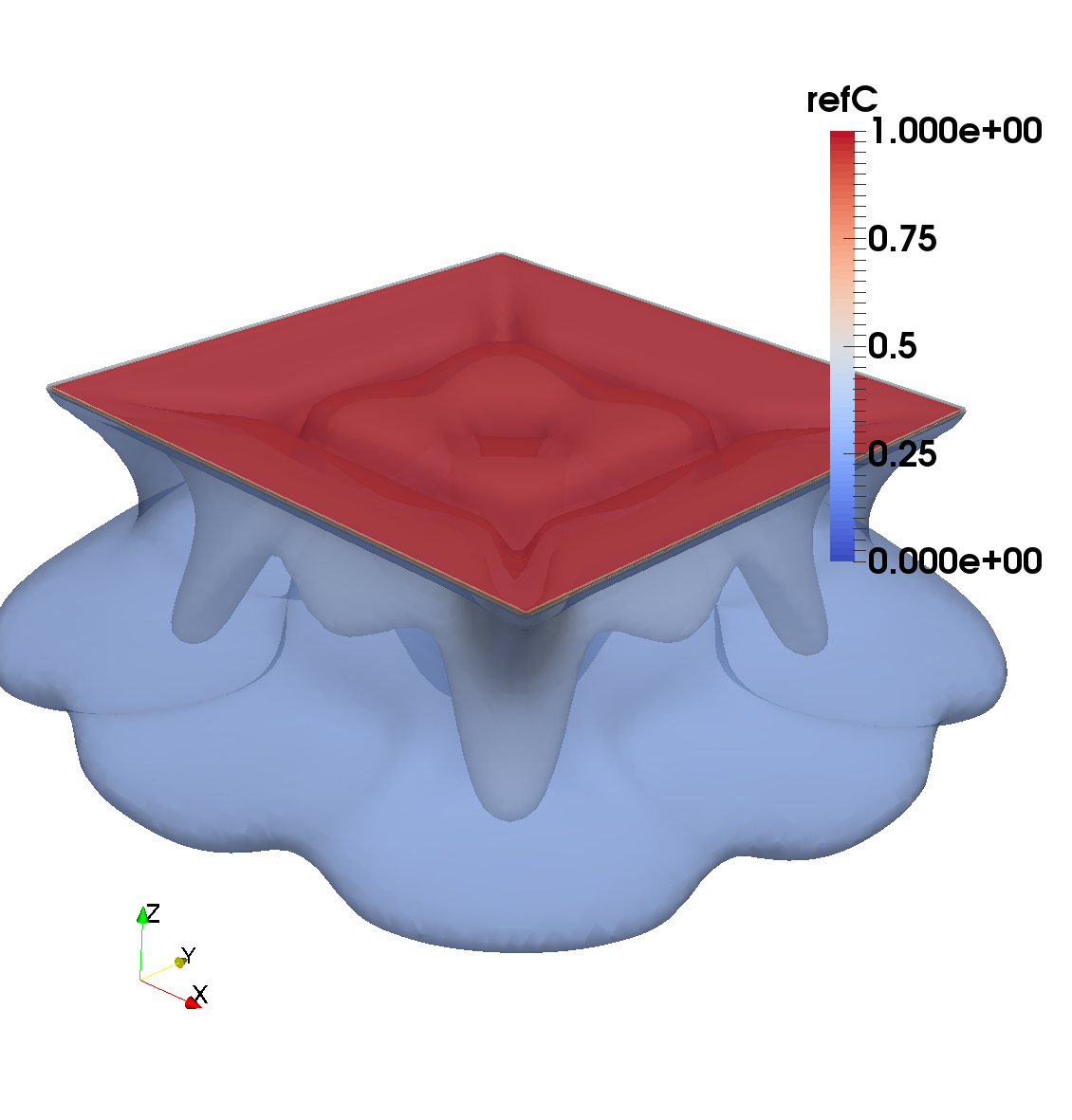}
     \label{subfig:refsol-up}
    \end{subfigure}
    \begin{subfigure}[b]{0.48\textwidth}
    \centering
    \caption{}
    \includegraphics[width=0.99\textwidth]{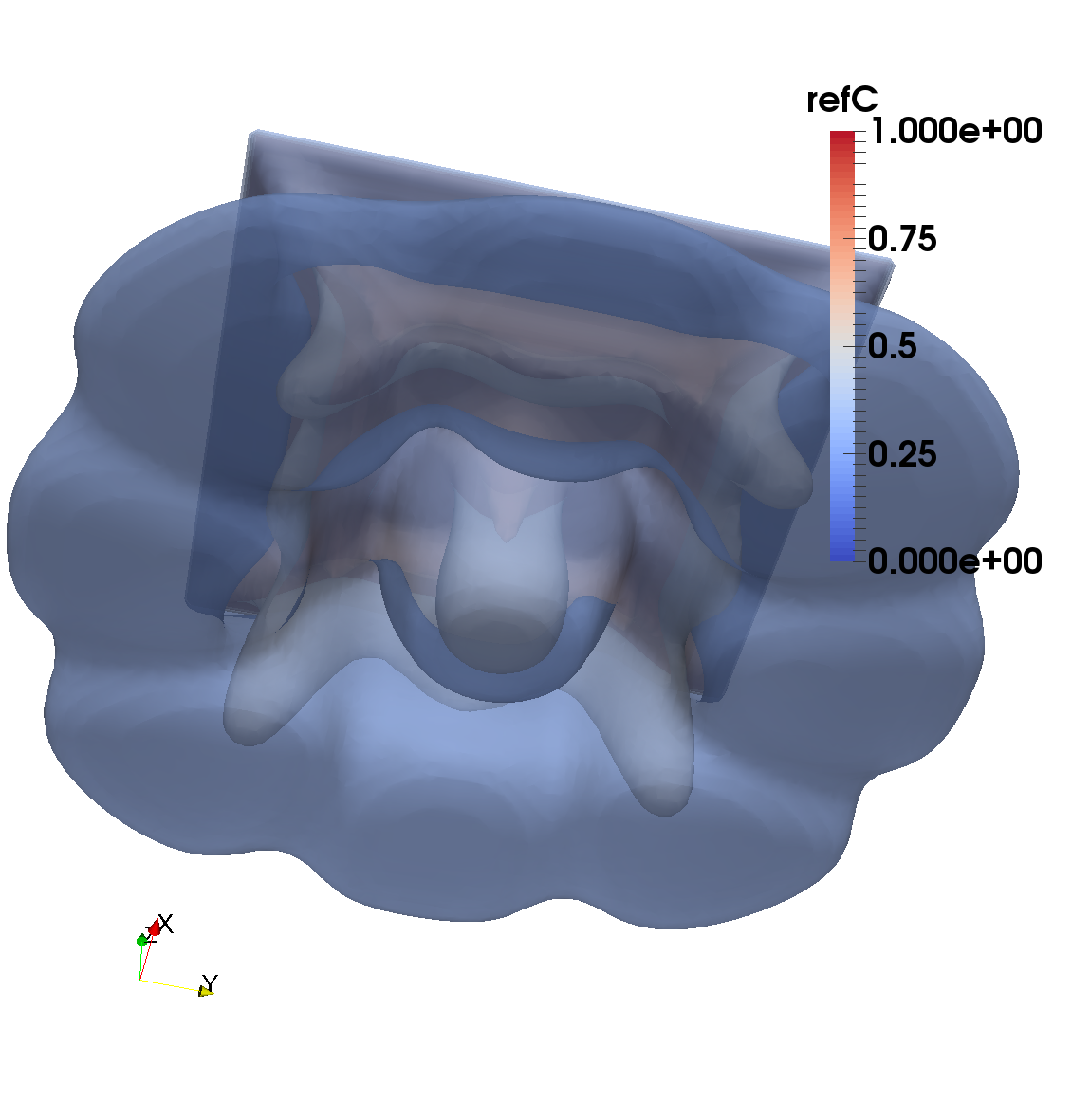}
     \label{subfig:refsol-bottom}
    \end{subfigure}
\caption{Five isosurfaces of the deterministic solution $\{\bx:\;{\sol}(t,\bx,0)\}=\{1.0,0.75,0.5,0.25,0.0\}$ after $t=9.6$ years; (a) view from the top; (b) view from the bottom.  The upper square edge corresponds to the boundary of the red area in Fig.\ \ref{fig:Elder3d-scheme}(a).}
\label{fig:3D_refsol}
\end{figure}

Five different contour lines (isosurfaces) of the deterministic solution $\sol(\bx,\btheta=0)$ for the porosity, defined in \refeq{eq:poro-paral3RVS}, are shown in Figures~\ref{fig:3D_refsol}. % and \ref{fig:domain3D_ts600_L5_M3_var}.
\newpage
\subsection{Elliptical cylinder with 3 RVs}
\label{ssec:cyl}
In this example we consider a cylindrical domain (see Fig.~\ref{subfig:Elder3d-cyl}). 
The goal is to see the influence of the computing domain on the solution. As we can see, the fingers have different shapes and numbers.

\begin{figure}[t]
\begin{center}
    \includegraphics[width=0.29\textwidth]{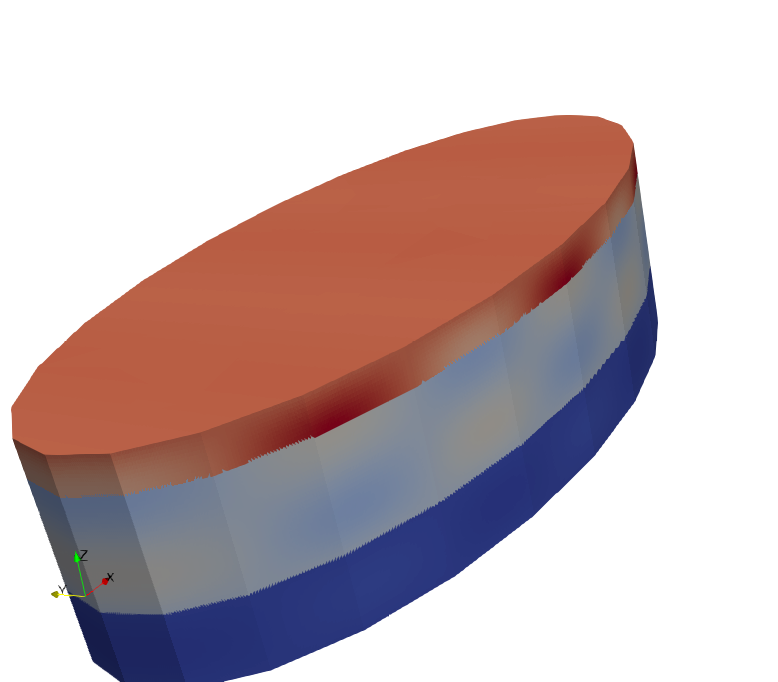}
    \includegraphics[width=0.29\textwidth]{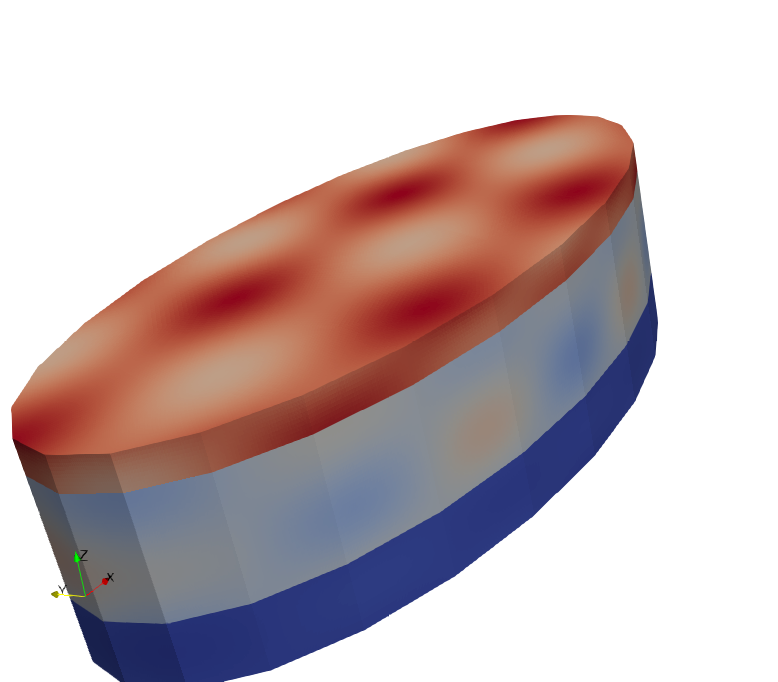}
    \caption{Two realizations of the porosity field, $\poro(\bx)\in [0.079,0.13]$}
\label{fig:poro_Older3d_26July_m3_cyl}
\end{center}
\end{figure}

The porosity coefficient has three layers and is defined with three RVs as follows:
\begin{align}
&\poro(t,\bx,\thetab)=0.1+0.05\cdot c_0\cdot \left( \frac{\xi_1 x}{600}\cos\frac{\pi x}{300} + \xi_2\sin\frac{\pi y}{150} + \xi_3\cos\frac{\pi x}{300}\sin\frac{\pi y}{150} \right) \label{eq:cyl_3layers} \\
 &\Bigg\{
 \begin{array}{lll}
c_0=0.01 &\quad \text{if} &z \leq -100\\ %+10\theta_1\\
c_0=0.10 &\quad \text{if}& -100 <z \leq -50\\% +10\theta_2\\
c_0=1.0 &\quad \text{if}& -50 <z \leq 0% +10\theta_2\\
\end{array}
\label{eq:cyl_3layers_}
\end{align}
Figure~\ref{fig:poro_Older3d_26July_m3_cyl} shows two random realizations of the porosity field.

\begin{figure}[t]
   \begin{subfigure}[b]{0.45\textwidth}
    \centering
     \caption{}
    \includegraphics[width=0.99\textwidth]{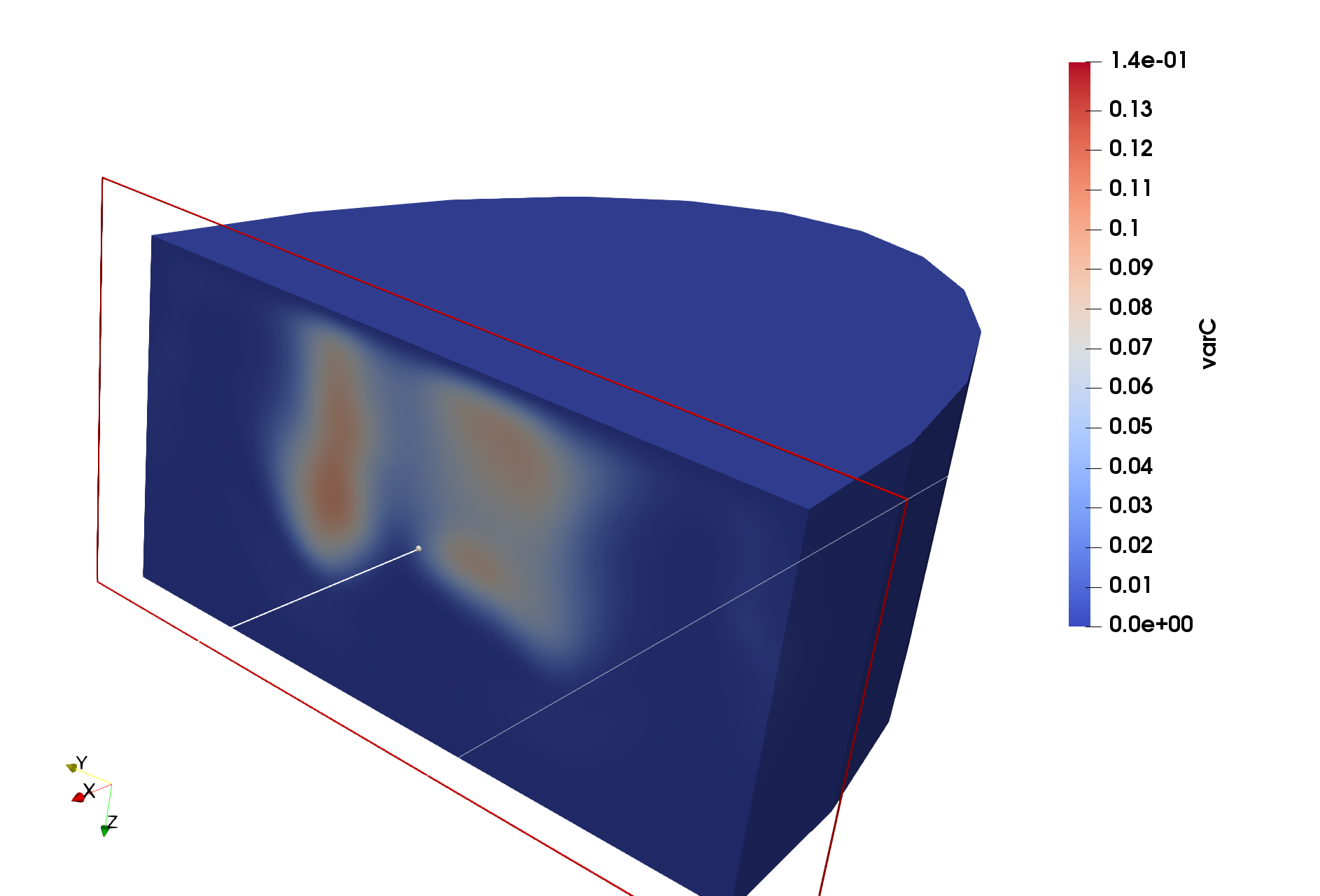}
     \label{subfig:MC200_var_cut1_ts400}
    \end{subfigure}
   \begin{subfigure}[b]{0.45\textwidth}
    \centering
     \caption{}
    \includegraphics[width=0.99\textwidth]{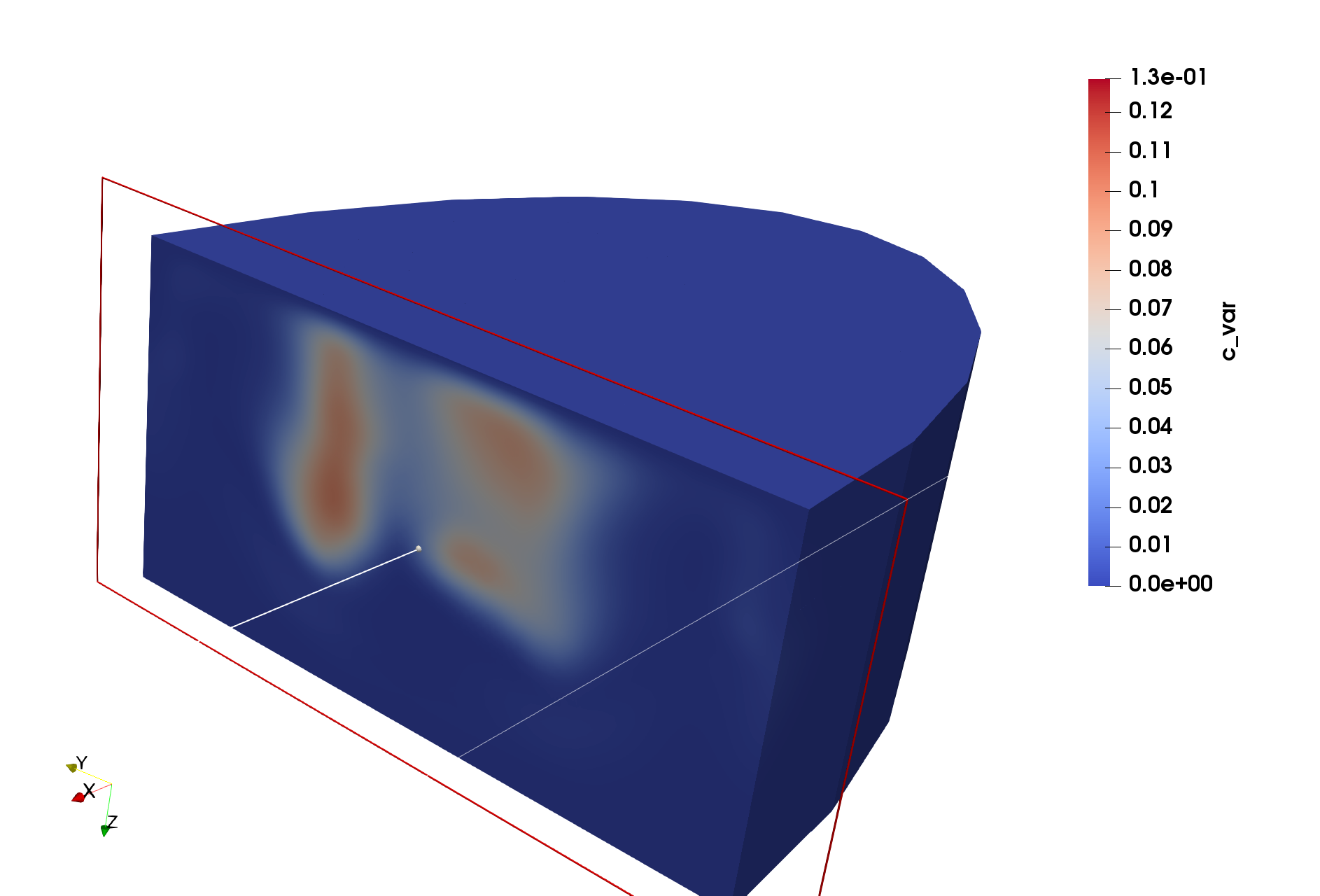}
     \label{subfig:PC39_var_cut1_ts400}
    \end{subfigure}\\
   \begin{subfigure}[b]{0.45\textwidth}
    \centering
     \caption{}
    \includegraphics[width=0.99\textwidth]{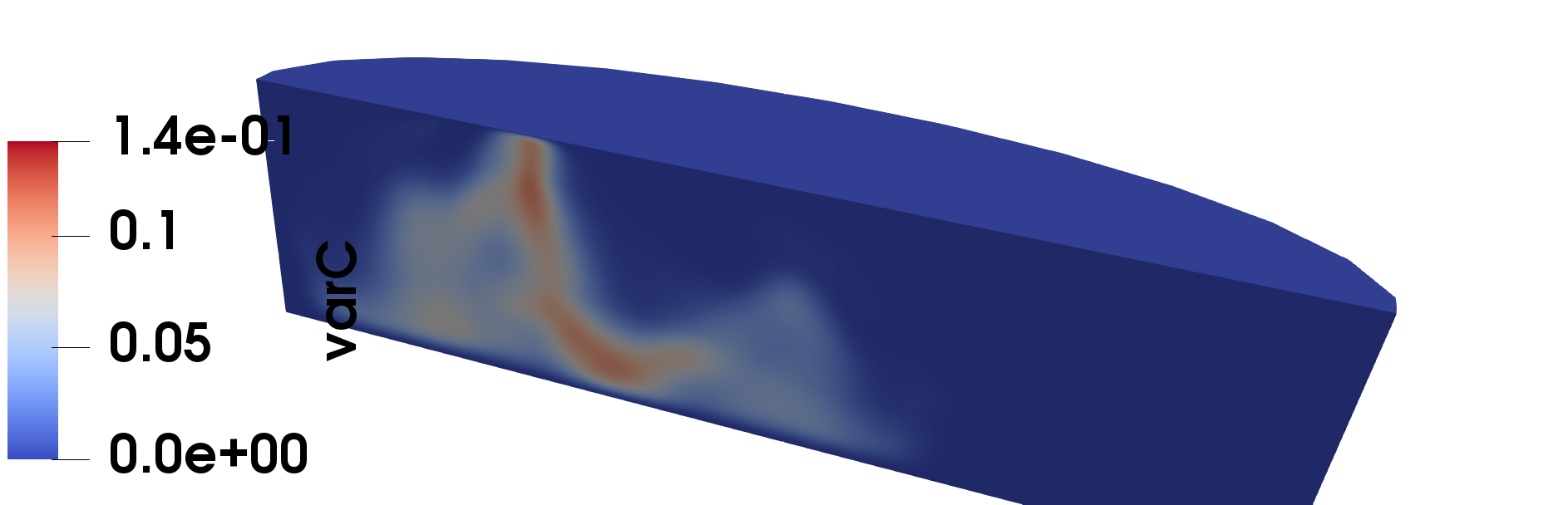}
     \label{subfig:MC200_var_cut3_ts400}
    \end{subfigure}
 %  \begin{subfigure}[b]{0.17\textwidth}
 %   \centering
 %   \includegraphics[width=0.99\textwidth]{figs4/video_cyl125_T400_3.png}
 %    \caption{}
 %    \label{subfig:v3}
 %   \end{subfigure}
   \begin{subfigure}[b]{0.45\textwidth}
    \centering
     \caption{}
    \includegraphics[width=0.99\textwidth]{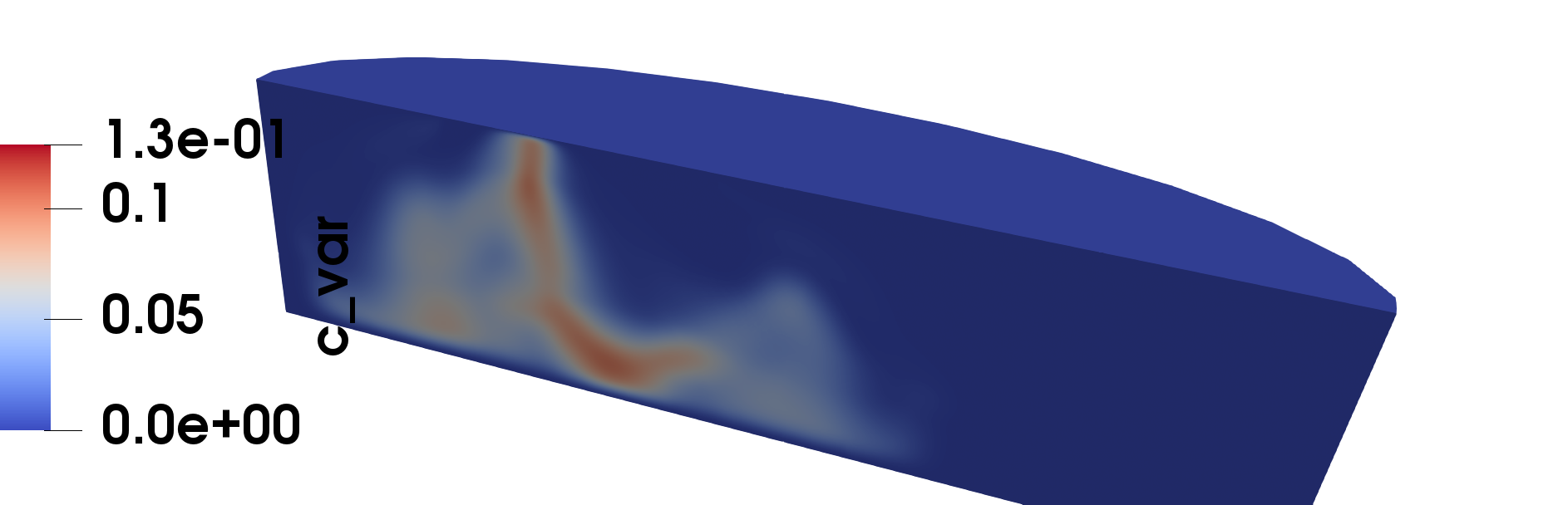}
     \label{subfig:PC35_var_cut3_TS400}
    \end{subfigure}
    \caption{Variance of the mass fraction in the cylindrical 3D reservoir, computed via 200 \QMC simulations (a) and (c) and \gPC of order 4 with 35 \gPC coefficients (b) and (d). $\var{\sol}_{\QMC}\in [0,0.14]$, $\var{\sol}_{\gPC}\in [0,0.13]$. The first cutting plane on (a)-(b) has normal vector $(1,0,0)$, the second on (c)-(d) $(0,1,0)$.}
    \label{fig:var_cut}
\end{figure}

In Fig.~\ref{fig:var_cut} we present a comparison of the variance of $\sol$, computed via \QMC (200 simulations) and \gPC of order $p=4$. All $\frac{(m+p)!}{m!p!}=\frac{7!}{3!4!}$=35 \gPC coefficients are computed with a full Clenshaw-Curtis quadrature, containing 125 quadrature points. In Figures~\ref{subfig:MC200_var_cut1_ts400} and \ref{subfig:PC39_var_cut1_ts400} the cutting plane has the normal vector $(1,0,0)$, and in Fig.~\ref{subfig:MC200_var_cut3_ts400} and \ref{subfig:PC35_var_cut3_TS400} vector $(0,1,0)$.
The variances, computed via \QMC and \gPC lie in the intervals $\var{\sol}_{\QMC}\in [0,0.14]$, and $\var{\sol}_{\gPC}\in [0,0.13]$ respectively. The small difference is not really visible in the pictures. Thus, we conclude that the \gPC surrogate model, computed on the full Clenshaw-Curtis quadrature grid can be used for approximating the mean and the variance of $\sol$.

\begin{figure}[t]
   \begin{subfigure}[b]{0.32\textwidth}
    \centering
     \caption{}
        \includegraphics[width=0.99\textwidth]{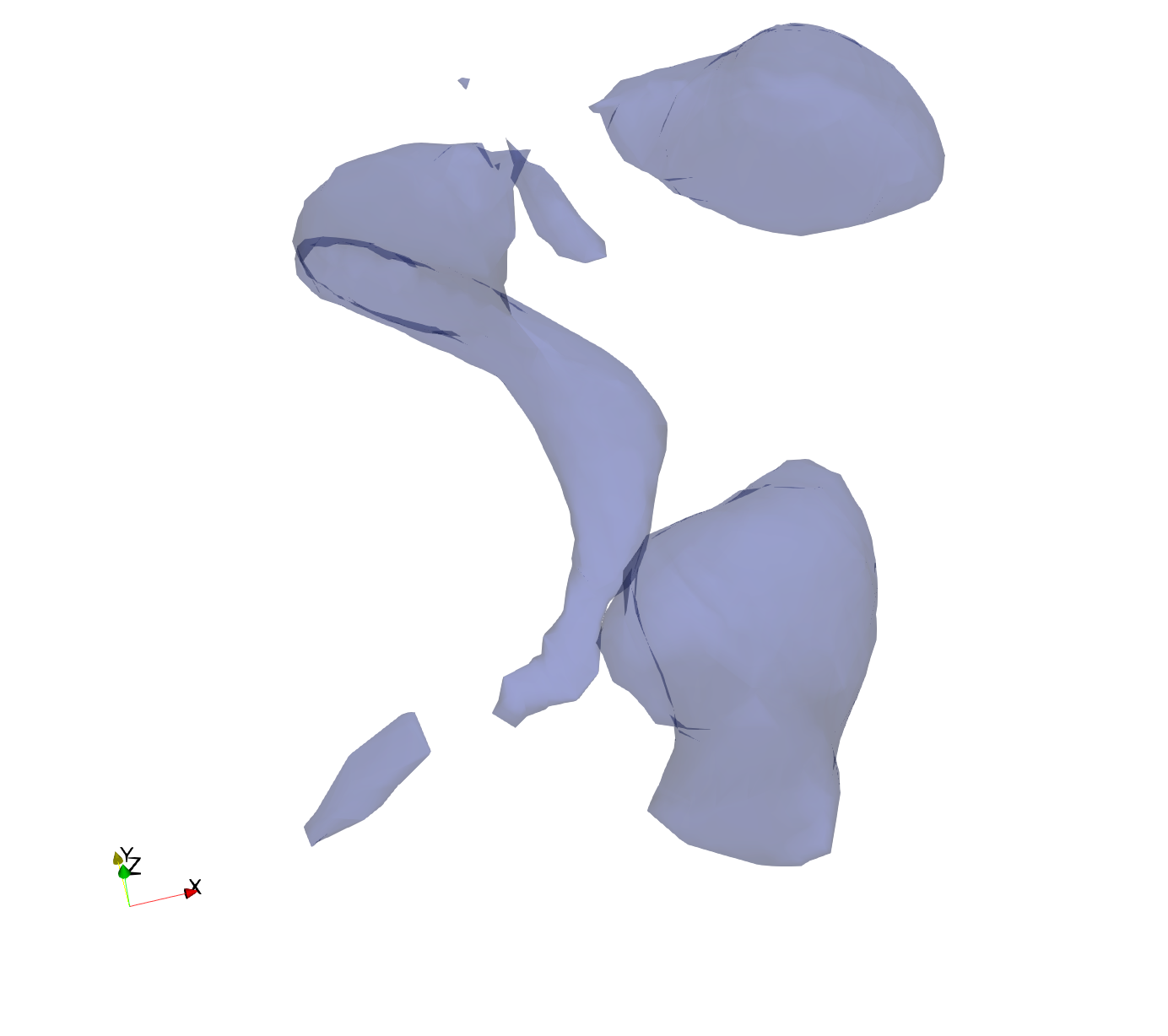}
     \label{subfig:var_c0_05_qmc}
    \end{subfigure}
   \begin{subfigure}[b]{0.32\textwidth}
    \centering
     \caption{}
        \includegraphics[width=0.99\textwidth]{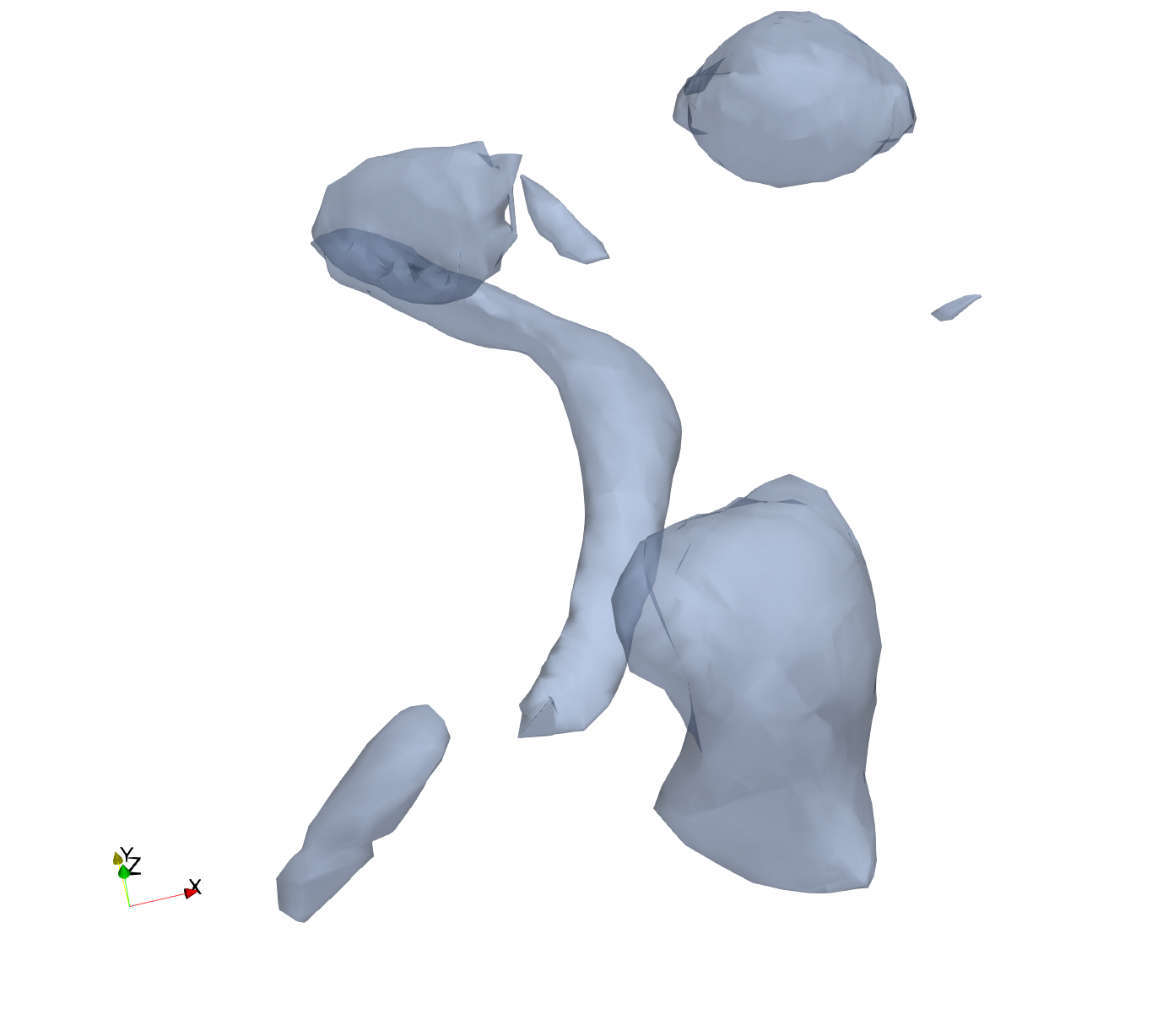}
 \label{subfig:var_c0_05_gpc}
    \end{subfigure}
   \begin{subfigure}[b]{0.32\textwidth}
    \centering
     \caption{}
    \includegraphics[width=0.99\textwidth]{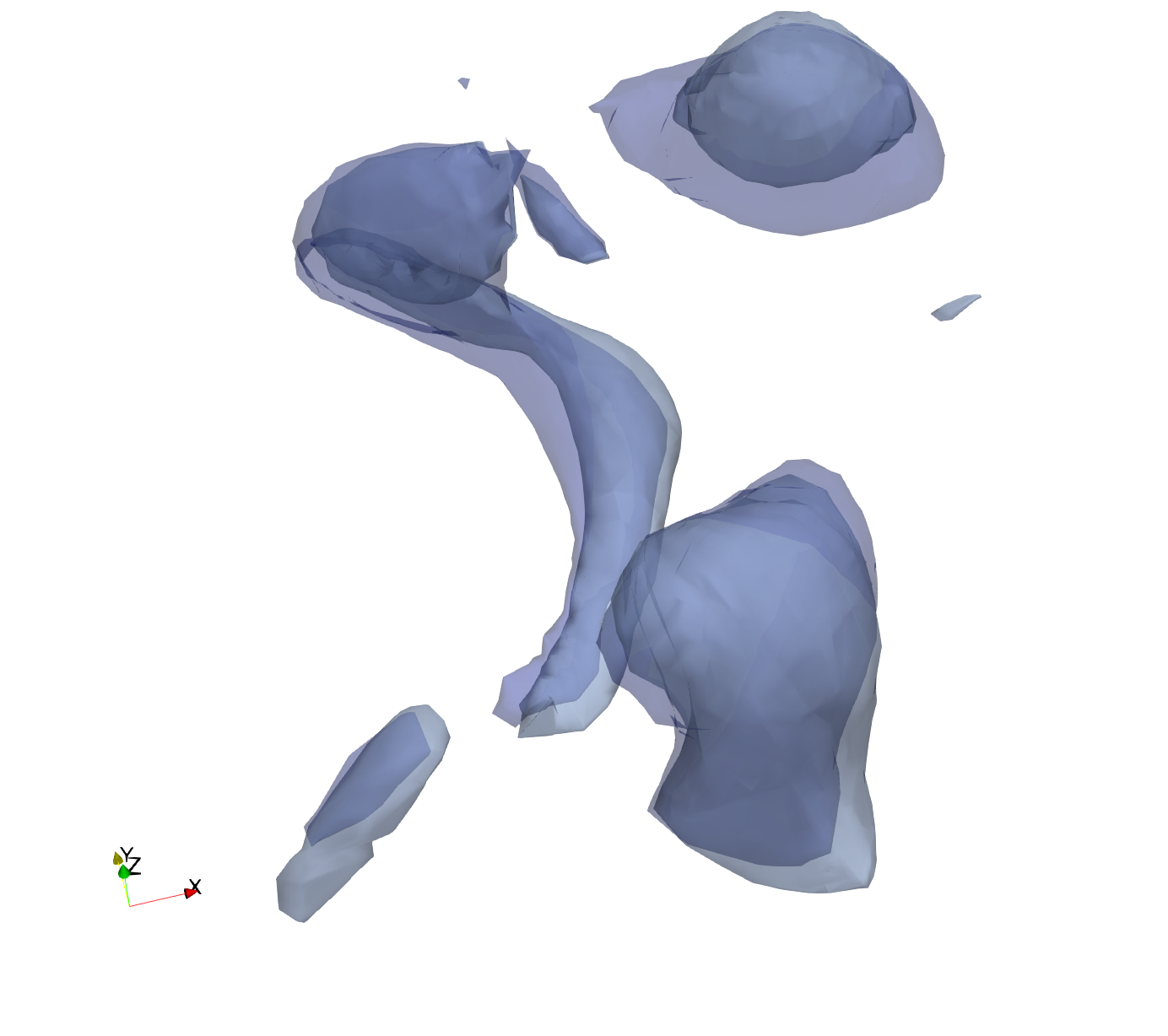}
     \label{subfig:var_cont_compare_0_05}
    \end{subfigure}
    \caption{The contour $\var{\sol}=0.025$, computed by a) \QMC, b) \gPC and c) comparison of both.}
    \label{fig:var_Older3d_26July_m3_cyl}
\end{figure}

\begin{figure}[htbp!]
\center
    \begin{subfigure}[b]{0.32\textwidth}
    \centering
      \caption{}
  \includegraphics[width=0.99\textwidth]{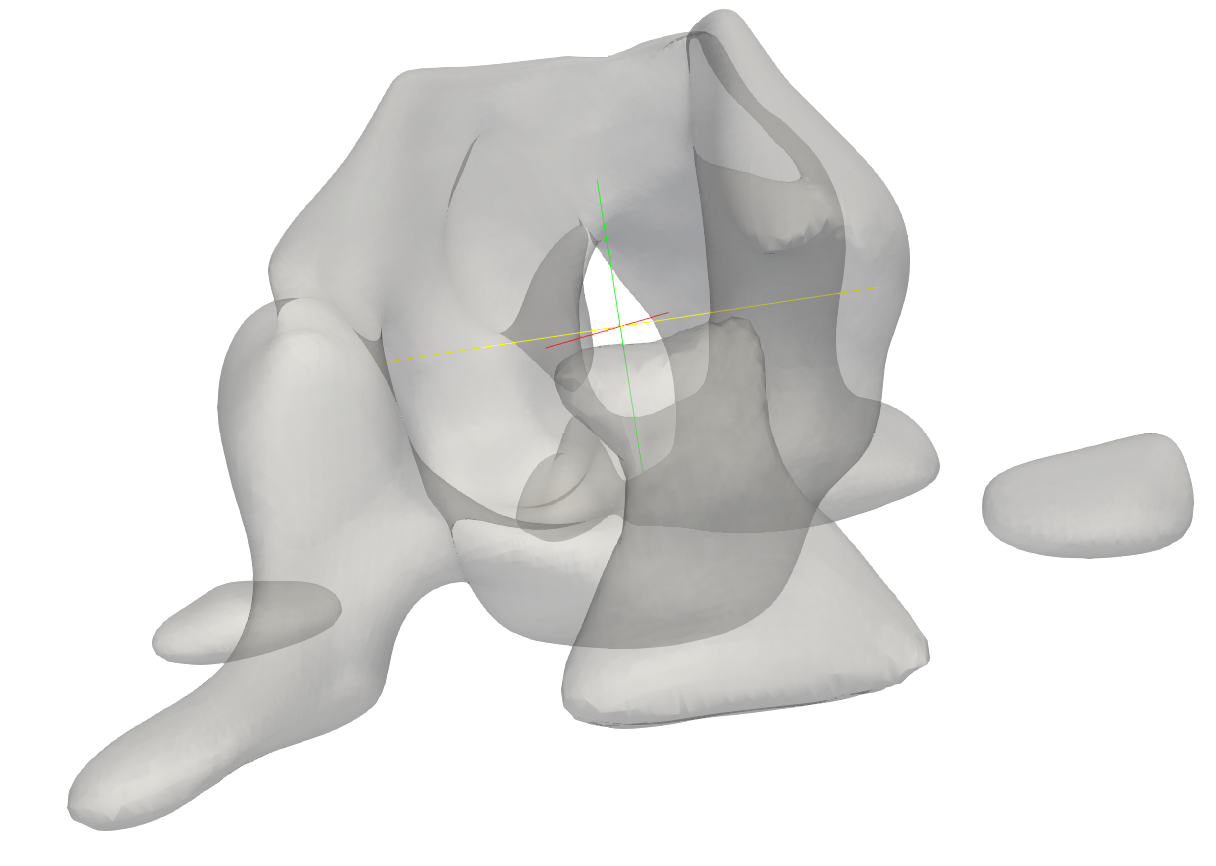}
     \label{subfig:QMC200_cont_0_07var}
    \end{subfigure}
  \begin{subfigure}[b]{0.32\textwidth}
    \centering
      \caption{}
  \includegraphics[width=0.99\textwidth]{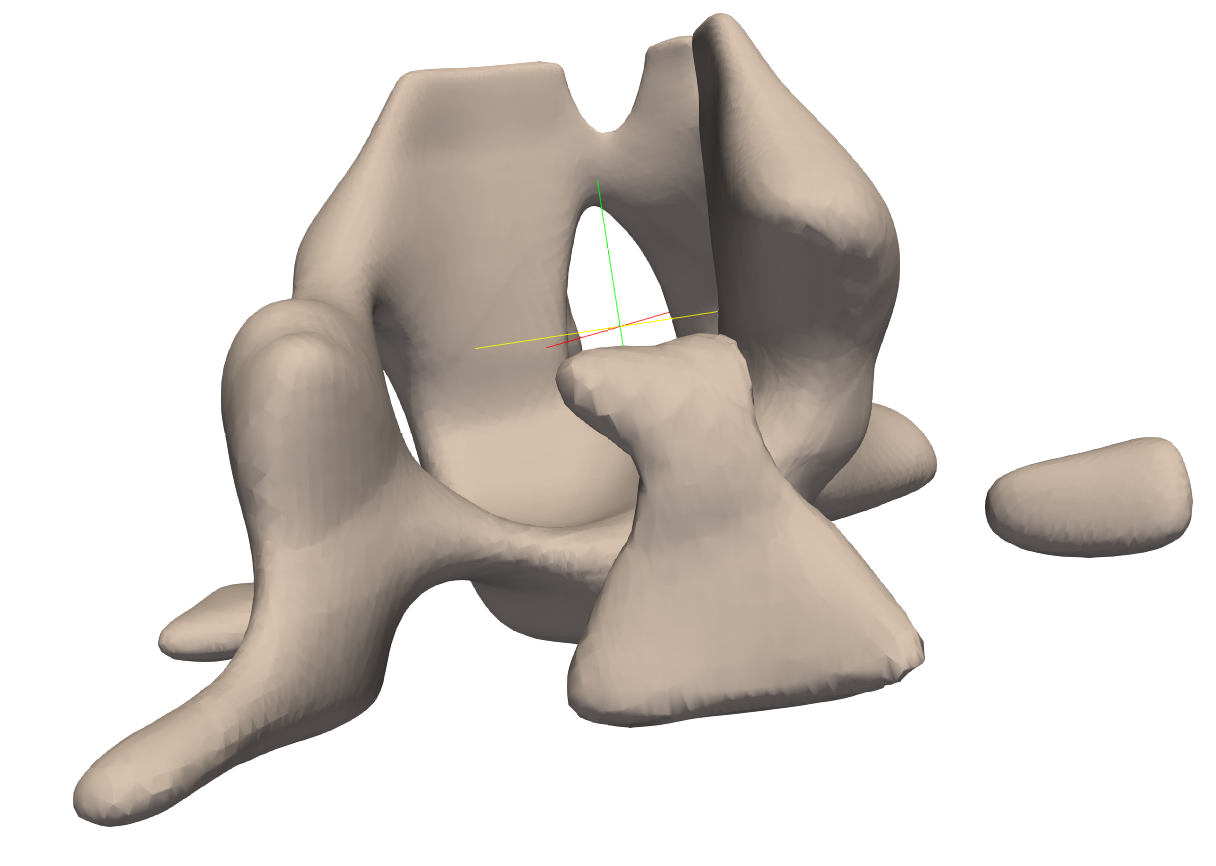}
     \label{subfig:PC4_35coefs_cont_0_07var}
    \end{subfigure}
  \begin{subfigure}[b]{0.32\textwidth}
    \centering
      \caption{}
  \includegraphics[width=0.99\textwidth]{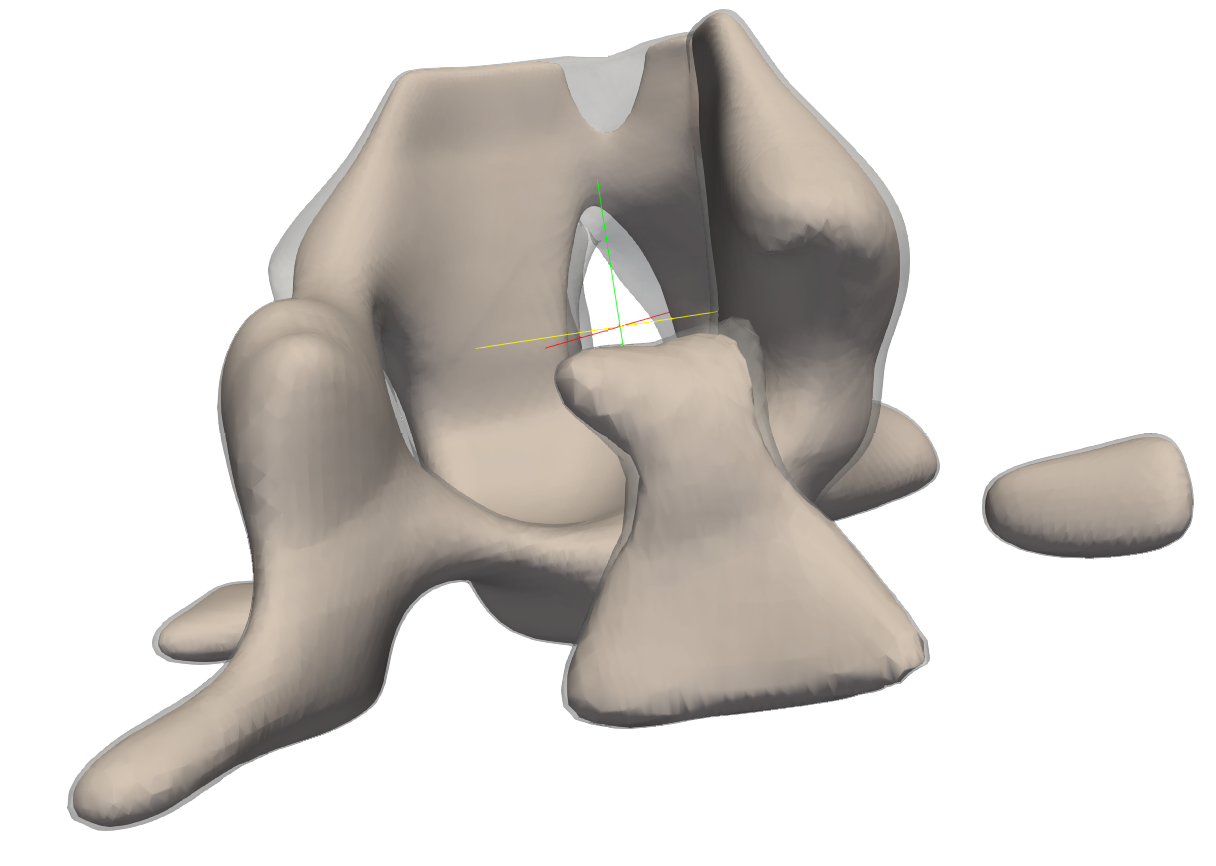}
     \label{subfig:Both_cont_0_07var}
    \end{subfigure}
\caption{Isosurface $\{\bx:\;\var{\sol}(\bx)=0.07\}$, computed via a) \QMC and b) \gPC response surface of order 4, and c) comparison of both isosurfaces}
\label{fig:3D_var_countur007}
\end{figure}

To demonstrate more similarities between the solutions, computed via \QMC and \gPC we visualize
isosurfaces of the variance of the mass fraction in the cylindrical 3D reservoir (Figures~\ref{fig:var_Older3d_26July_m3_cyl} and \ref{fig:3D_var_countur007}). 
% The isosurface on Fig.~\ref{fig:3D_var_countur007} is computed via 200 \QMC simulations. 
% The isosurface on Fig.~\ref{subfig:isoPC} is computed via \gPC of order 4 with 35 \gPC coefficients. 
% We also have $\var{\sol}_{\QMC}\in [0,0.14]$ and $\var{\sol}_{\gPC}\in [0,0.13]$.
% \begin{figure}[t]
%   \begin{subfigure}[b]{0.45\textwidth}
%     \centering
%      \caption{}
%     \includegraphics[width=0.99\textwidth]{figs/Older_cyl_MC_var_0_14.png}
%      \label{subfig:isoMC}
%     \end{subfigure}
%   \begin{subfigure}[b]{0.45\textwidth}
%     \centering
%      \caption{}
%     \includegraphics[width=0.99\textwidth]{figs/Older_cyl_PC_var_0_13.png}
%      \label{subfig:isoPC}
%     \end{subfigure}
%     \caption{Isosurfaces of the variation of the mass fraction in the cylindrical 3D reservoir. a) $\sol$ is computed via 200 \QMC simulations and (b) $\sol$ is computed via \gPC of order 4 with 35 \gPC coefficients. $\var{\sol}_{\QMC}\in [0,0.14]$, $\var{\sol}_{\gPC}\in [0,0.13]$.}
%     \label{fig:var_iso}
% \end{figure}
\newpage
In Figure~\ref{fig:var_Older3d_26July_m3_cyl} we demonstrate two isosurfaces $\{\bx:\;\var{\sol}(\bx)=0.025\}$ computed via a) \QMC method, and b) \gPC of order 4 surrogate. The Figure~\ref{subfig:var_cont_compare_0_05} shows the comparison of both isosurfaces. As we can see the results are very similar, the isosurface computed via \QMC is slightly larger.

In Figure~\ref{fig:3D_var_countur007} we demonstrate two isosurfaces $\{\bx:\;\var{\sol}(\bx)=0.07\}$ computed via a) \QMC method, Fig.~\ref{subfig:QMC200_cont_0_07var}, and b) \gPC of order 4 surrogate, Fig.~\ref{subfig:PC4_35coefs_cont_0_07var}. The Figure~\ref{subfig:Both_cont_0_07var} shows the comparison of both isosurfaces. As we can see the results are very similar, the isosurface computed via \QMC is slightly larger.

\begin{figure}[htbp!]
\center
    \begin{subfigure}[b]{0.32\textwidth}
    \centering
      \caption{}
  \includegraphics[width=0.99\textwidth]{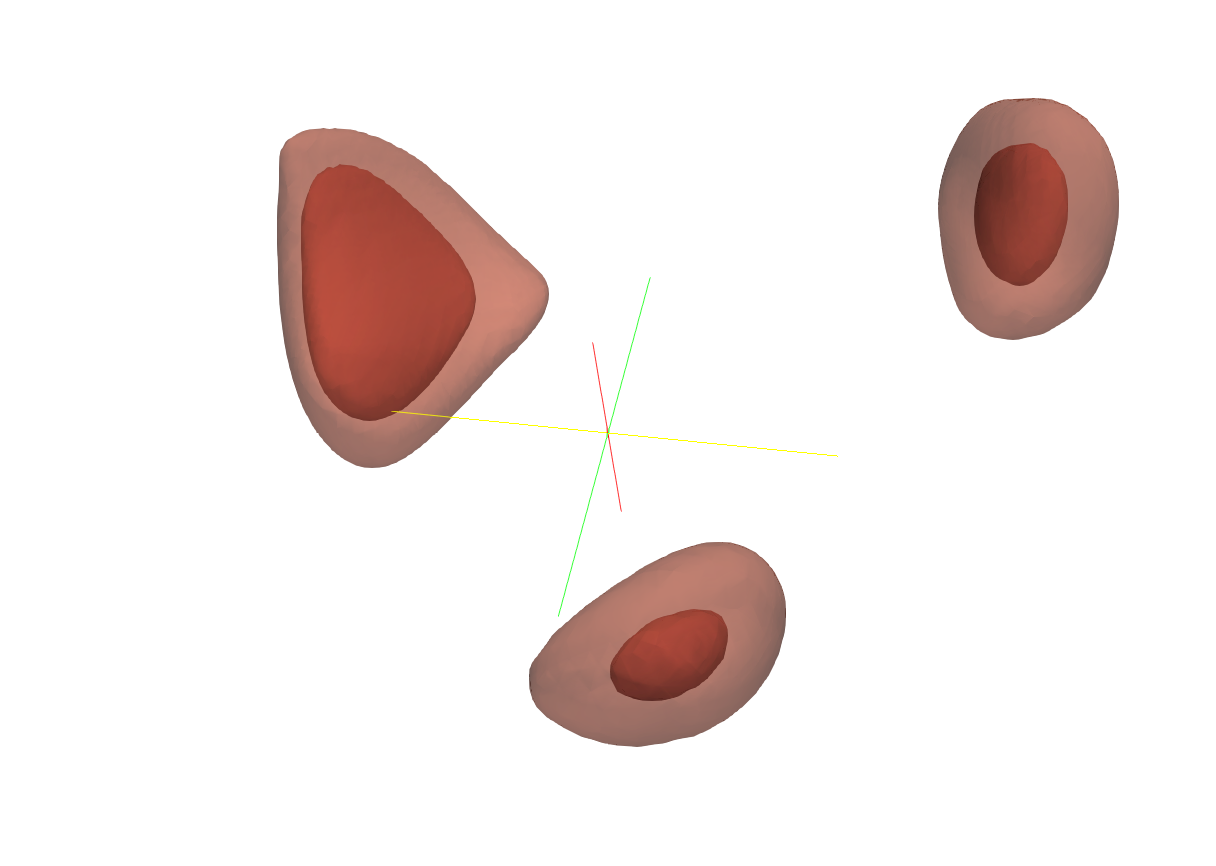}
     \label{subfig:cont_compar_var0_12}
    \end{subfigure}
\caption{Comparison of isosurfaces $\{\bx:\;\var{\sol}(\bx)=0.12\}$, computed via \QMC and \gPC response surface of order 4.
The isosurface computed via \QMC is slightly larger.}
\label{fig:Both_cont_0_12}
\end{figure}

In Figure~\ref{fig:Both_cont_0_12} we compare two isosurfaces $\{\bx:\;\var{\sol}(\bx)=0.12\}$ computed via \QMC method and \gPC of order 4 surrogate. As we can see the results (shape and number) are very similar again. The isosurface computed via \QMC is slightly larger.

%
%
% Figure~\ref{fig:var_3layers_Iso} shows two isosurfaces: $\var{\sol}_{0.05}$ and $\var{\sol}_{0.15}$ after 3 years.
% \begin{figure}[htbp!]
%     \begin{subfigure}[b]{0.48\textwidth}
%     \centering
%     \caption{}
%     \includegraphics[width=0.99\textwidth]{figs/16July/var_gPCE_m3_p3_Iso0_05.png}
%      \label{subfig:var_gPCE_m3_Iso005}
%     \end{subfigure}
%     \begin{subfigure}[b]{0.48\textwidth}
%     \centering
%     \caption{}
%     \includegraphics[width=0.99\textwidth]{figs/16July/var_gPCE_m3_p3_Iso0_15.png}
%      \label{subfig:var_gPCE_m3_Iso015}
%     \end{subfigure}    
% \caption{Isosurfaces of the variance of the mass fraction after 3 years; (a) $\var{\sol}_{0.05}$; (b)  $\var{\sol}_{0.15}$.} 
% \label{fig:var_3layers_Iso}
% \end{figure}
%
\newpage
\newpage
\subsection{Time evolution of $\overline{\sol}(t,\bx)$ and $\var{\sol}(t,\bx)$}
\begin{figure}[t]
   \begin{subfigure}[b]{0.45\textwidth}
    \centering
     \caption{}
    \includegraphics[width=0.99\textwidth]{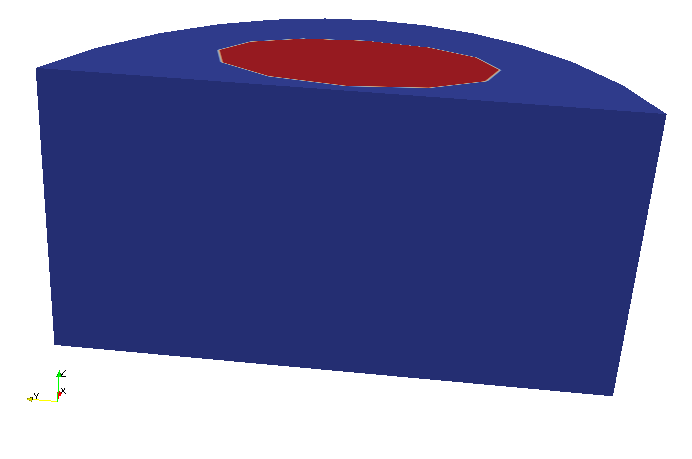}
     \label{subfig:v0}
    \end{subfigure}
   \begin{subfigure}[b]{0.45\textwidth}
    \centering
     \caption{}
    \includegraphics[width=0.99\textwidth]{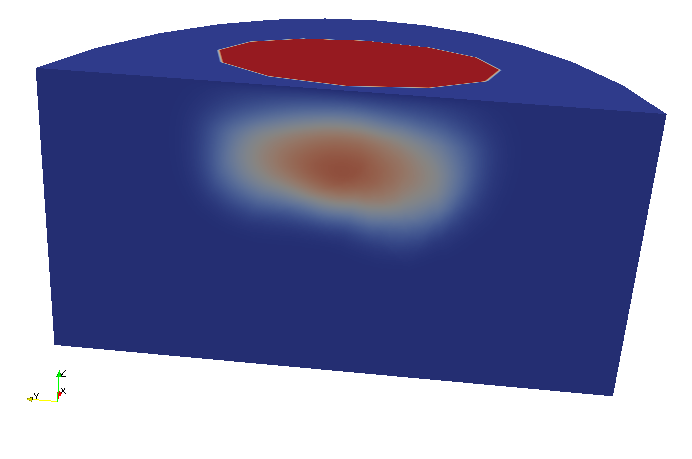}
     \label{subfig:v1}
    \end{subfigure}\\
   \begin{subfigure}[b]{0.45\textwidth}
    \centering
     \caption{}
    \includegraphics[width=0.99\textwidth]{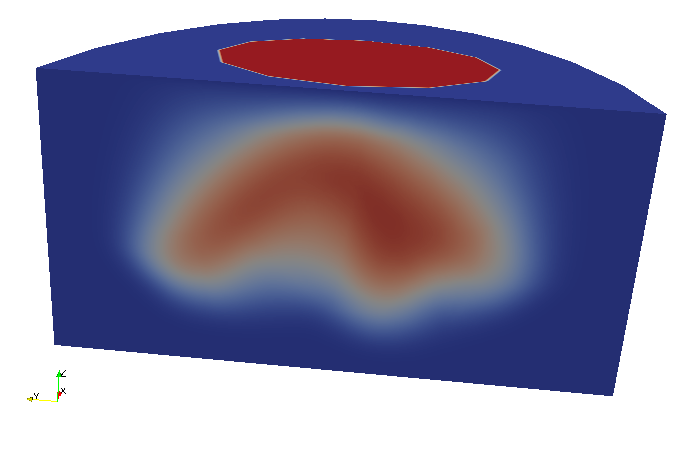}
     \label{subfig:v2}
    \end{subfigure}
 %  \begin{subfigure}[b]{0.17\textwidth}
 %   \centering
 %   \includegraphics[width=0.99\textwidth]{figs4/video_cyl125_T400_3.png}
 %    \caption{}
 %    \label{subfig:v3}
 %   \end{subfigure}
   \begin{subfigure}[b]{0.45\textwidth}
    \centering
     \caption{}
    \includegraphics[width=0.99\textwidth]{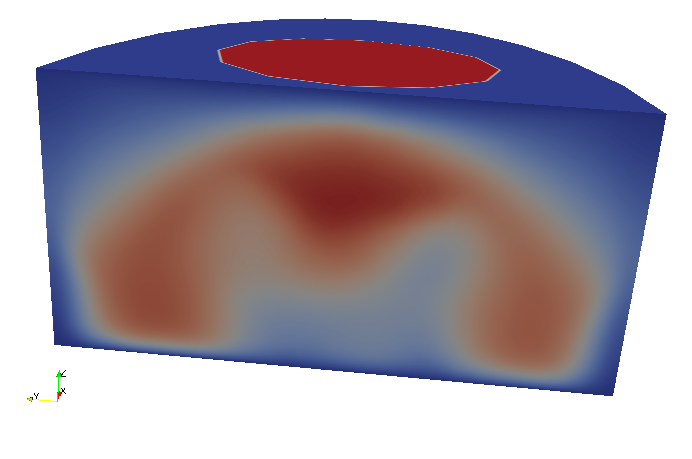}
     \label{subfig:v4}
    \end{subfigure}
    \caption{Propagation of the mean of the mass fraction in the cylindrical 3D reservoir, cutting plane is $(150,y,z)$. Evolution of $\EXP{\sol}$ in time after 
  a) 0, b) 0.55, c) 1.1, d) 2.2 years.}
    \label{fig:video_Older3d}
\end{figure}

\begin{figure}[t]
   \begin{subfigure}[b]{0.32\textwidth}
    \centering
     \caption{}
    \includegraphics[width=0.99\textwidth]{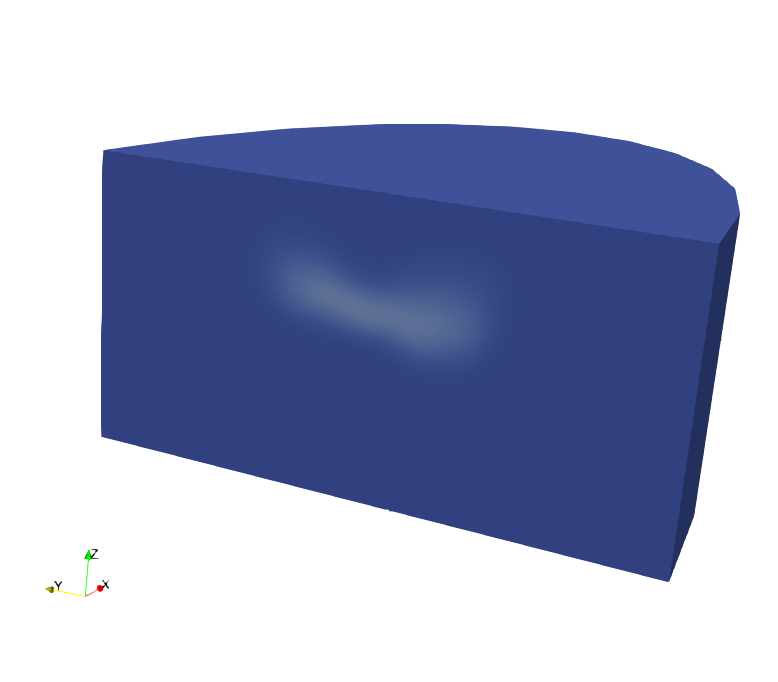}
     \label{subfig:v0_t100}
    \end{subfigure}
   \begin{subfigure}[b]{0.32\textwidth}
    \centering
     \caption{}
    \includegraphics[width=0.99\textwidth]{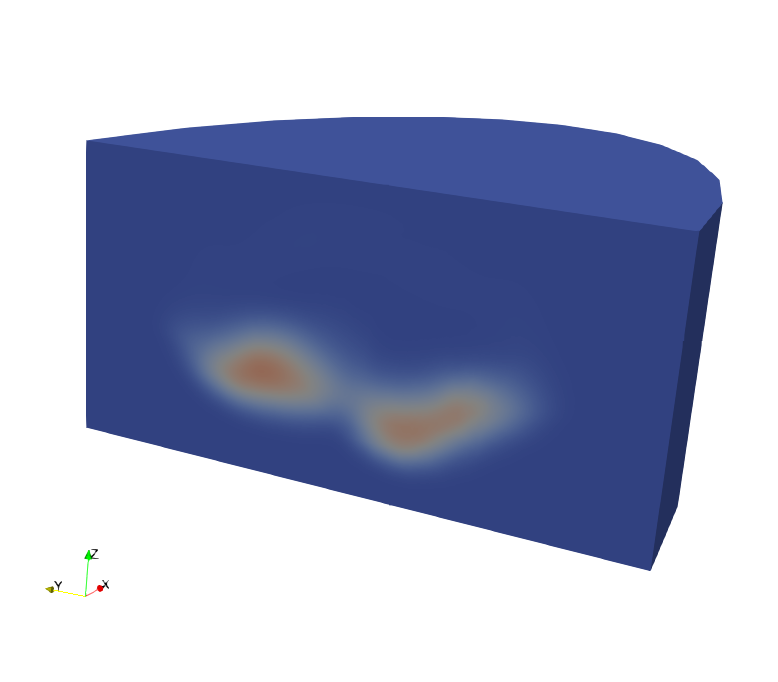}
     \label{subfig:v1_t200}
    \end{subfigure}
%   \begin{subfigure}[b]{0.32\textwidth}
%     \centering
%      \caption{}
%     \includegraphics[width=0.99\textwidth]{figs/cyl/var_time_evol/var_T300_cyl_MC.png}
%      \label{subfig:v2_t300}
%     \end{subfigure}
   \begin{subfigure}[b]{0.32\textwidth}
    \centering
     \caption{}
    \includegraphics[width=0.99\textwidth]{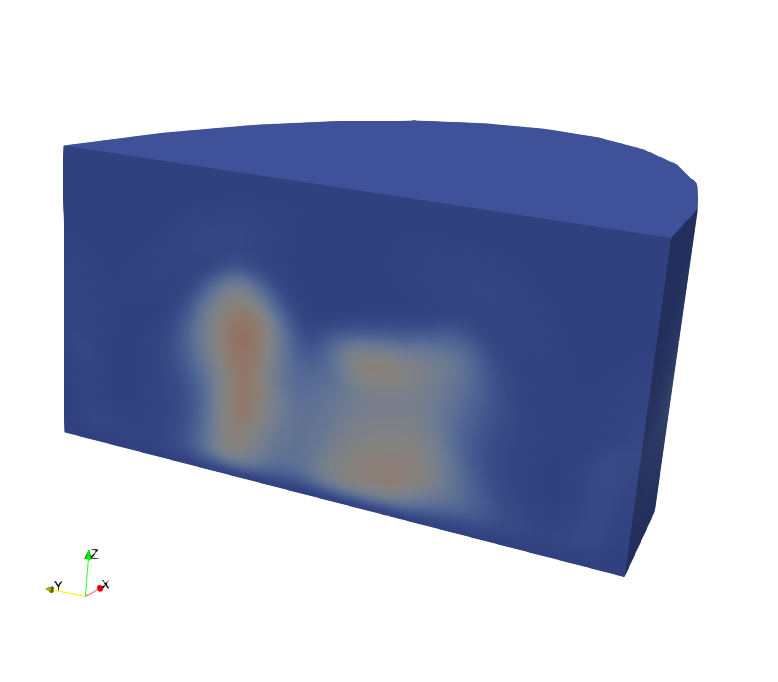}
     \label{subfig:v4_t400}
    \end{subfigure}
    \caption{Propagation of the variance of the mass fraction in the cylindrical 3D reservoir, cutting plane is $(150,y,z)$. Evolution of the mean concentration in time after 
    a) 0.55, b) 1.1, d) 2.2 years (100, 200, 400 time steps respectively). $\var{\sol}\in[0,0.17]$.}
    %$\{a) 0, b) 0.55, c) 1.1, 1.65, 2.2\}$ years.}
    \label{fig:video_Older3d_var}
\end{figure}

This is again an experiment with $\poro$ as in \refeq{eq:cyl_3layers}-\refeq{eq:cyl_3layers_}. The \gPC coefficients were computed with 125 Gauss-Legendre quadrature points.
Figure~\ref{fig:video_Older3d} shows evolution of the mean of the mass fraction in the cylindrical 3D reservoir after a) 0, b) 0.55, c) 1.1, d) 2.2 years. The cutting plane is $(150,y,z)$. One can clearly see different stages of the ''finger" building.
Figure~\ref{fig:video_Older3d_var} shows 
the variance of the mass fraction in the cylindrical 3D reservoir after $\{0.55, 1.1,  2.2\}$ years (after 100, 200, and 400 time steps, respectively).

\subsection{Three layers example with cylindrical domain}
The porosity coefficient inside of each level is random and is defined by the following equation
\begin{minipage}{0.6\textwidth}
\begin{align*}
    \poro(t, \bx,\btheta)&=0.1+c_z(\theta_1\cos(\pi x/600)+\theta_2\cos(\pi y/300)\\
    &+\theta_3\sin(\pi x/600)\cos(\pi z/150)),
\end{align*}    
with $c_z=\bigg\{
 \begin{array}{cc}
0.01, & -150 \leq z < -100\\
0.1, &  -100 \leq z < -50\\
1.0, &  -50 \leq z \leq 0.\\
\end{array}
$

Here, $\btheta=(\theta_1,\theta_2,\theta_3)$ is a quadrature point on a full Gaussian-Legendre grid (level 3, dimension 3).
\end{minipage}
\begin{minipage}{0.39\textwidth}
\begin{center}
    \quad \quad \includegraphics[width=0.69\textwidth]{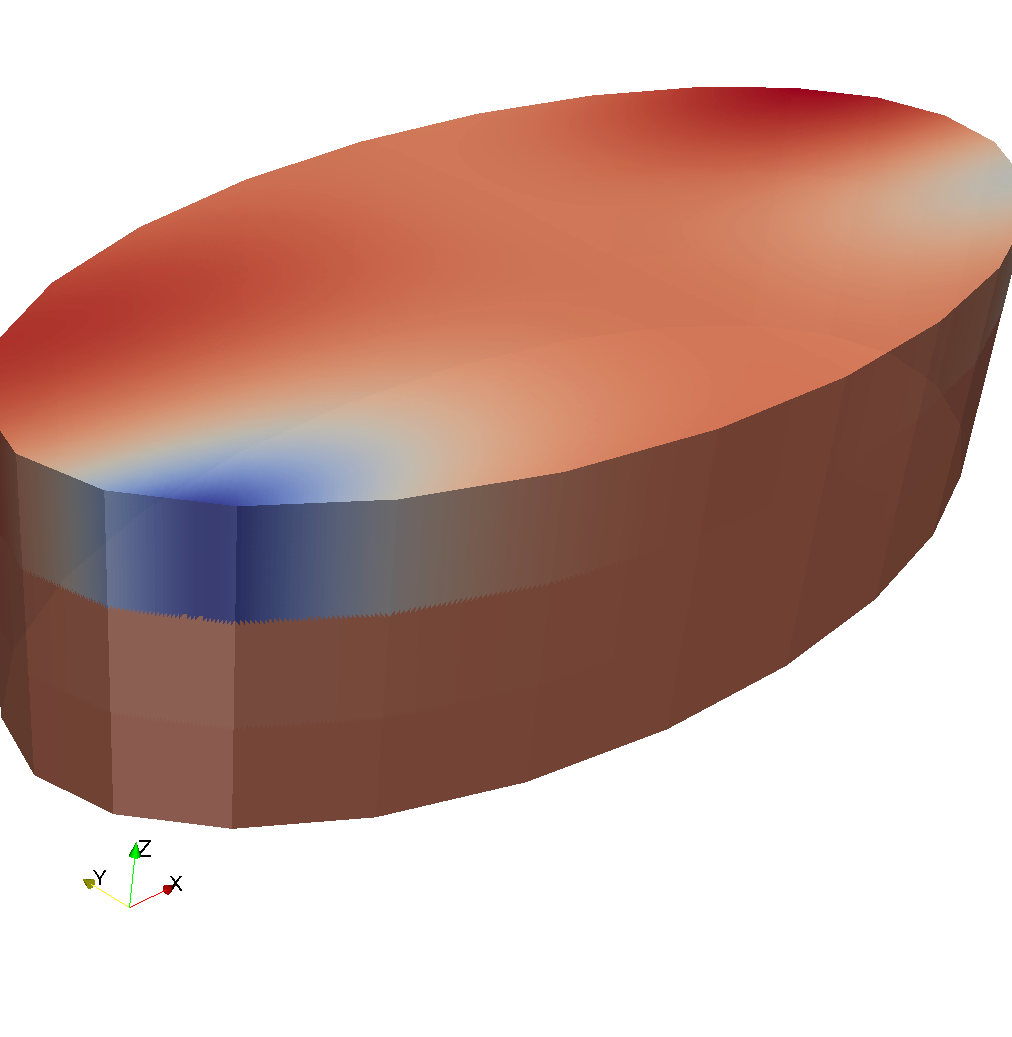}\\
    Fig. One realization of the porosity field, $\poro\in[0.05,0.13]$.
    \end{center}
    \vspace{0.5cm}
\end{minipage}

The maximum and the minimum over all porosity samples and over all $\bx$ are $0.02$ and $0.15$, respectively. One random realization of the porosity is shown in Figure on the right.
The random porosity has three horizontal layers along $z$-axes: the first $z\in [-150,-100)$, second $z\in [-100,-50)$, and third $z\in[-50,0]$.

\begin{figure}[htbp!]
    \begin{subfigure}[b]{0.33\textwidth}
    \centering
    \includegraphics[width=0.99\textwidth]{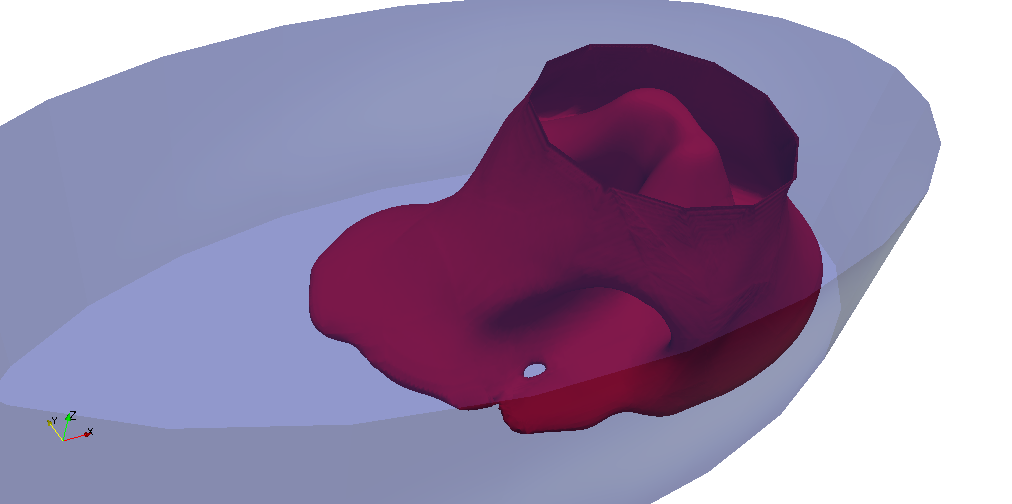}
     \caption{}
     \label{subfig:mean_gPCE_Iso065}
    \end{subfigure}
    \begin{subfigure}[b]{0.33\textwidth}
    \centering
    \includegraphics[width=0.99\textwidth]{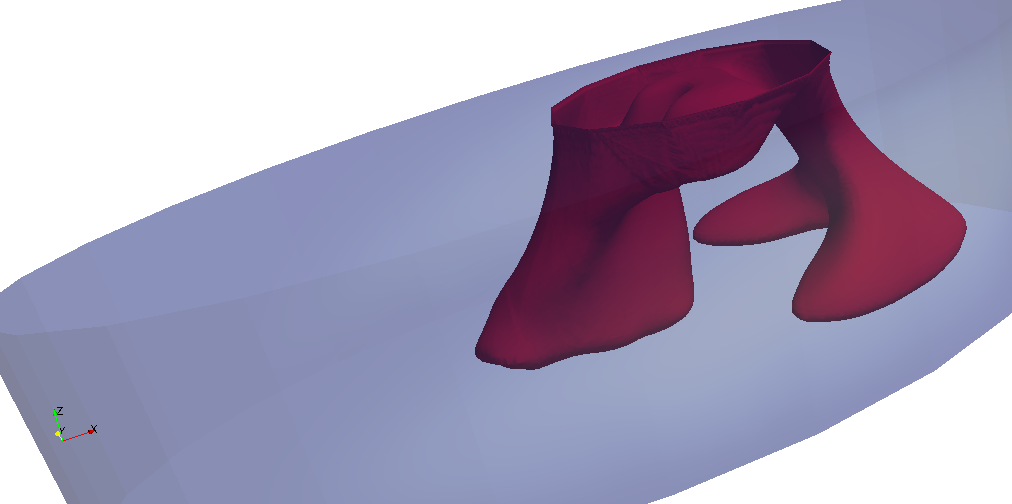}
     \caption{}
     \label{subfig:mean_gPCE_Iso078}
    \end{subfigure}    
    \begin{subfigure}[b]{0.33\textwidth}
    \centering
    \includegraphics[width=0.99\textwidth]{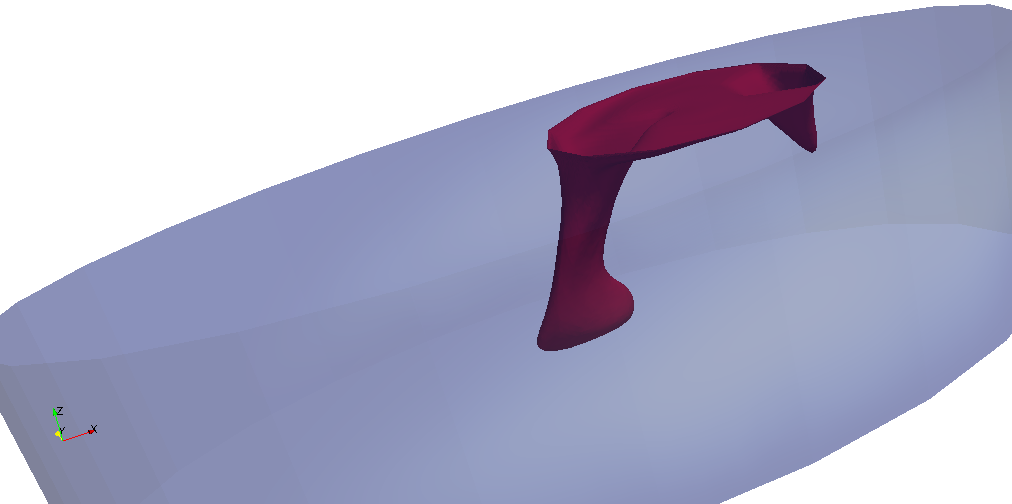}
     \caption{}
     \label{subfig:mean_gPCE_Iso095}
    \end{subfigure}    
\caption{Isosurfaces of the mean of the mass fraction after 500 time steps ($\approx 2.5$ years); (a) $\{\bx:\;\overline{\sol}(\bx)=0.65\}$, (b) $\{\bx:\;\overline{\sol}(\bx)=0.78\}$, (c) $\{\bx:\;\overline{\sol}(\bx)=0.98\}$.}
\label{fig:mean_3layers}
\end{figure}

Figure~\ref{fig:mean_3layers} shows three different isosurfaces of the mean value of the mass fraction: (a) $\{\bx:\;\overline{\sol}(\bx)=0.65\}$, (b) $\{\bx:\;\overline{\sol}(\bx)=0.78\}$, and (c) $\{\bx:\;\overline{\sol}(\bx)=0.95\}$. The large blue area visualizes the reservoir contour and the smaller dark isosurface represents a isosurface of the mean value.

\begin{figure}[htbp!]
    \begin{subfigure}[b]{0.33\textwidth}
    \centering
    \includegraphics[width=0.99\textwidth]{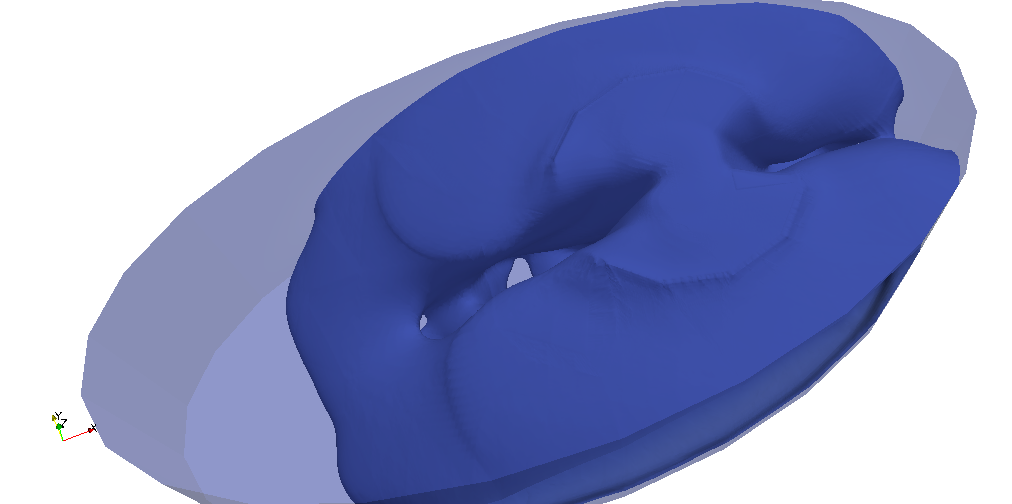}
     \caption{}
     \label{subfig:ar_gPCE_m3_p3_Iso0_01}
    \end{subfigure}    
    \begin{subfigure}[b]{0.33\textwidth}
    \centering
    \includegraphics[width=0.99\textwidth]{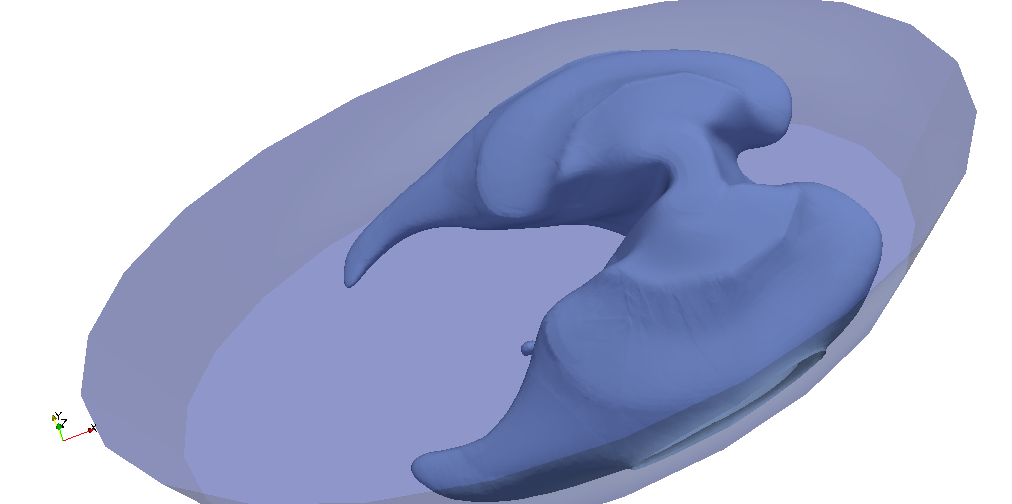}
     \caption{}
     \label{subfig:ar_gPCE_m3_p3_Iso0_05}
    \end{subfigure}    
    \begin{subfigure}[b]{0.33\textwidth}
    \centering
    \includegraphics[width=0.99\textwidth]{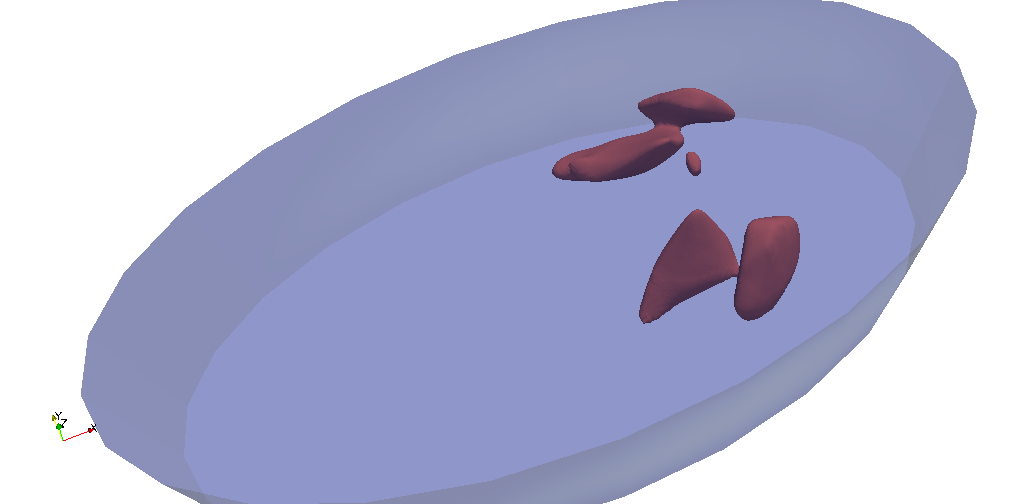}
     \caption{}
     \label{subfig:var_gPCE_m3_p3_Iso0_15}
    \end{subfigure}    
\caption{Isosurfaces of the variance of the mass fraction after 500 time steps, $\{\bx:\;\var{\sol}(\bx)=0.01\}$, (b) $\{\bx:\;\var{\sol}(\bx)=0.05\}$, (c) $\{\bx:\;\var{\sol}(\bx)=0.15\}$.} 
\label{fig:var_3layers3}
\end{figure}

Three isosurfaces of the variance of the mass fraction after 500 time steps ($\approx 2.5$ years) are shown in Fig.~\ref{fig:var_3layers3}. These isosurfaces are $\{\bx:\;\var{\sol}(\bx)=0.01\}$, (b) $\{\bx:\;\var{\sol}(\bx)=0.05\}$, and (c) $\{\bx:\;\var{\sol}(\bx)=0.15\}$, respectively. The \gPC surrogate model of order $p=3$ with $M=3$ random variables was used to compute the mean and the variance of the mass fraction.

\newpage
\subsection{Best practices}
\label{sec:best}
After numerous numerical tests, we created the best practice guide as below:
\begin{enumerate}
\item The number of uncertain variables to describe the porosity coefficient depends on the model and assumptions. We recommend to start with 1-3 variables to understand the phenomena better. 
We note that a heterogeneous porosity field in 3D may require hundreds of random variables. That means thousands of \gPC coefficients. This is computationally un-affordable unless one uses some advance low-rank tensor techniques \cite{dolgov2014computation} or multi-scale/homogenization strategy. In this work, we considered three random variables. The more random variables are present, the harder it could be to interpret the results, and additional sensitivity analysis may be required.
\item It is assumed a smooth dependence of the output QoI on the input uncertain parameters. We recommend to check it at least visually before computing \gPC.
\item The \gPC-based surrogate, for the cases when the variance of the porosity is not so high, showed reliable results. For the validation, one can use a \QMC method.
\item We were not able to get a correct variance of $\sol$ with the \gPC method in tests where the variance of the porosity was large.
\item We recommend the maximal polynomial order in \gPC to be $p=2,3,4$ or $5$. With our computational budget, we failed to get the correct variance of $\sol$ for $p>5$. Probably the quadrature rule was not good enough to catch oscillatory behavior of \gPC polynomials. This is a general known drawback of using global polynomials.
\item The coefficients of the \gPC surrogate require calculating high-dimensional integrals. For stochastic dimensions 4 and higher, we recommend to use Gauss-Legendre (or Clenshaw-Curtis) sparse grid. For dimensions 1-3, we recommend to use a full grid. 
\item The polynomial order in \gPC is 2-4. Otherwise, the number of simulations will be too large.
\item The total numerical cost could be reduced by choosing the spatial and time step sizes adaptive, depending on the random perturbation. But in our implementation, it is hard to it do on the fly. Therefore, we recommend ``conservative'' settings, which are the same for all random simulations. A large time or large spatial step sizes may result in lack of convergence.   
\end{enumerate}
\newpage
\section{Conclusion}
\label{litv:sec:Conclusion}

We solved the density driven groundwater flow problem with uncertain porosity and permeability. An accurate solution of this time-dependent and non-linear problem is impossible due to the presence of natural uncertainties in the porosity and permeability in the reservoir. 
Therefore, we estimated the mean value and the variance of the solution, as well as the propagation of uncertainties from the random input parameters to the solution. For this, we, first, computed a multi-variate polynomial approximations (\gPC) to the solution and then used it to estimate the required statistics. 
Utilizing the \gPC method allows us 
to reduce the computational cost compare of the classical Monte Carlo method.
\gPC assumes that the output function $\sol(t,\bx,\thetab)$ is square-integrable and smooth w.r.t uncertain input variables $\btheta$.

We considered two different aquifers - a solid parallelepiped and a solid elliptic cylinder. One of our goals was to see how the domain geometry influences the shape of fingers.

Since the considered problem is nonlinear, 
a high variance in the porosity may result in totally different solutions, for instance, the number of fingers, their intensity, the shape, the propagation time, and the velocity may strongly variate.

From the implementation point of view, we have demonstrated how to use the parallel multigrid solver \myug as a black-box solver.
An additional novelty of this work is that all algorithms for uncertainty quantification are implemented on a distributed memory system, where all solutions and \gPC coefficients are distributed over multiple nodes. We saw great potential in using the \myug library as a black box solver for quantification of uncertainties. The results are reproducible, and the code is available online on GitHub repository\footnote{\url{https://github.com/UG4},\quad \url{https://github.com/UG4/ughub}}.  

The number of cells in the presented experiments varied from $241{,}152$ till $15{,}433{,}728$ for the cylindrical domain and from $524{,}288$ till $4{,}194{,}304$ for parallelepiped. The maximal number of parallel processing units was $800\times 32$, where $800$ is the number of parallel nodes, and $32$ is the number of computing cores on each node. The computing time varied from 2 hours for the coarse mesh to 23 hours for the finest mesh (for one simulation).

Although this work let many questions open, it is an important step towards the investigation of the density driven flow phenomena.
Potential further applications are an estimation of the contamination risks, monitoring the quality of drinking water, modeling of the seawater intrusion into coastal aquifers, radioactive waste disposal, and contaminant plumes.

As a next step, we will integrate the \myug multi-grid solver in the Multi Level Monte Carlo framework with the idea to compute mostly samples on coarse grids and just a few on a fine grid. 

Another idea is to infer the (unknown) porosity coefficient. We assume that there are some noisy measurements of $\sol$ at some locations $\{\bx_1,\bx_2,\ldots,\bx_k\}$ are available. Then the prior probability density function of $\sol$ could be improved via Bayesian update formula \cite{
hgmEzBvrAlOp2016,hgmEzBvrAl:2016,Rosic2013,bvrAlOpHgm11,litvinenko2013inverse}. Herewith, the solution $\sol(\bx^*)$ at a point $\bx^*$ (or at $k$ points or in a sub-domain) could be computed with a much smaller numerical cost and with the smaller memory requirement as the whole solution in the whole domain \cite{Litvinenko17_PatrInv,LitDiss}.
% Note, that these statistics could be further used for more efficient Bayesian inference, data assimilation, the optimal design of experiment, and optimal control. The overall goal is modeling of water pollution and monitoring of the quality of subsurface water flow.%

\textbf{Acknowledgements:}\\
We would like to thank KAUST HPC support team for the assistance with Shaheen II.
This work was supported by the Extreme Computing Research Center, by the SRI-UQ Strategic Initiative and by Computational Bayesian group at King Abdullah University of Science and Technology.
%
%

%\bibliographystyle{siam}
% WHERE IS THIS FILE ? \bibliographystyle{spmpsci}
%\bibliography{new_article_Sydney_about_sampling, matthies_BU_paper-1, mybib, references}
%\bibliography{litvinenko_dolgov_khoromskij}
%\bibliography{MoCa}
\begin{footnotesize}

\end{footnotesize}

\begin{appendices}
\section{Legendre polynomials}
Since the uncertain input parameters have Uniform distribution, we employ multivariate Legendre polynomials as the \gPC basis.
Multivariate Legendre polynomials are defined on $[-1,1]^M$, where $M$ is the dimensionality of the stochastic space. The first few Legendre polynomials are 
\begin{align}
\label{eq:LegendreMonom}
\pol_0(x)&=1,\quad \pol_1(x)=x,\quad \pol_2(x)=\frac{3}{2}x^2-\frac{1}{2},\\
\pol_3(x)&=\frac{5}{2}x^2-\frac{3}{2}x,\quad 
\pol_4(x)=\frac{35}{8}x^4-\frac{30}{8}x^2+\frac{3}{8}\\
\pol_5(x)&=\frac{63}{8}x^5-\frac{70}{8}x^3+\frac{15}{8}x
\end{align}
The Legendre polynomials are orthogonal with respect to the $L_2$ norm on the interval $[-1,1]$.
The recursive formula is
\begin{equation}
\pol_0(x)=1, \quad \pol_1(x)=x,\quad (n+1)\pol_{n+1}(x)=(2n+1)x\pol_n(x) - n\pol_{n-1}(x).
\end{equation}
The scalar product is 
\begin{equation}
\int_{-1}^1 \pol_n(x)\pol_m(x)\mu(x)dx=\frac{1}{2n+1}\delta_{nm},
\end{equation}
or
\begin{equation}
\int_{-1}^1 \pol_n(x)\pol_m(x)dx=\frac{2}{2n+1}\delta_{nm},
\end{equation}
where the density function $\mu(x)$ is a constant (is equal 1/2 on $[-1,1]$), and where $\delta_{nm}$ denotes the Kronecker delta. 
The $L^2[-1,1]$ norm $\Vert \pol_n\Vert =\sqrt{\frac{2}{2n+1}}$ and the weighted $L^{2}_{\mu}[-1,1]$ norm is $\Vert \pol_n \Vert =\sqrt{\frac{1}{2n+1}}$. The \gPC coefficients can be computed as follows
\begin{equation}
c_{\alpha}(x):=\frac{(c(\thetab),\Pol_{\alpha}(\thetab))_{L^2_{\mu}}}{\Vert \Pol_{\alpha} \Vert^2_{L^2_{\mu}}}=
\frac{\int_{-1}^1 c\Pol_{\alpha} \mu(\thetab) d\thetab}{\int_{-1}^1 \Pol_{\alpha}\Pol_{\alpha} \mu(\thetab) d\thetab}=
\frac{\mu\int_{-1}^1 c\Pol_{\alpha} d\thetab}{\mu\int_{-1}^1 \Pol_{\alpha}\Pol_{\alpha} d\thetab}=\frac{0.5\cdot \sum_{i=1}^{n_q} c(\thetab_i)\Pol_{\alpha}(\thetab_i)w_i} {\frac{1}{2{\alpha}+1}}
\end{equation}
and for $d>1$ we have
\begin{equation}
c_{\alpha}(x):=\frac{(c(\thetab),\Pol_{\alpha}(\thetab))_{L^2_{\mu}}}{\Vert \Pol_{\alpha} \Vert^2_{L^2_{\mu}}}=
\frac{\int_{[-1,1]^d} c\Pol_{\alpha} \mu(\thetab) d\thetab}{\int_{[-1,1]^d} \Pol_{\alpha}\Pol_{\alpha} \mu(\thetab) d\thetab}=
\frac{\mu\int_{[-1,1]^d} c\Pol_{\alpha} d\thetab}{\mu\int_{[-1,1]^d} \Pol_{\alpha}\Pol_{\alpha} d\thetab}=\frac{0.5^d\cdot \sum_{i=1}^{n_q} c(\thetab_i)\Pol_{\alpha}(\thetab_i)w_i} {\prod_{i=1}^d\frac{1}{2\alpha_i+1}}.
\end{equation} 
\end{appendices}

\end{document}